\documentclass[leqno,12pt]{article}
\usepackage{amsmath,amssymb,a4wide}
\usepackage[curve,matrix,arrow,cmtip]{xy}

\usepackage{verbatim}

\NoComputerModernTips

\newdir^{ (}{{}*!/-3pt/\dir^{(}}    
\newdir_{ (}{{}*!/-3pt/\dir_{(}}    

\def\OM{\mathchoice
  {\rlap{\kern3.2pt$\overline{\phantom{L}}$}M}
  {\rlap{\kern3.2pt$\overline{\phantom{L}}$}M}
  {\rlap{\kern2.4pt$\scriptstyle\overline{\phantom{L}}$}M}
  {\rlap{\kern1.8pt$\scriptscriptstyle\overline{\phantom{L}}$}M}}

\let\le\leqslant
\let\ge\geqslant
\let\leq\leqslant

 \def\End{\mathop{\rm End}\nolimits}

 \def\Proj{\mathop{\rm Proj}\nolimits}
 \def\Spec{\mathop{\rm Spec}\nolimits}
 
 \def\deg{\mathop{\rm deg}\nolimits}
 
 \def\mod{\mathop{\rm mod}\nolimits}
\def\Ker{\mathop{\rm Ker}\nolimits}
\def\graph{\mathop{\rm graph}\nolimits}

\def\GL{\mathop{\rm GL}\nolimits}

\def\id{{\rm id}}

\let\phi\varphi
\let\epsilon\varepsilon
\let\setminus\smallsetminus

\newtheorem{Thm}{Theorem}[section]
\newtheorem{Prop}[Thm]{Proposition}
\newtheorem{Lem}[Thm]{Lemma}

\newtheorem{Def}[Thm]{Definition}

\newtheorem{Rem}[Thm]{Remark}

\def\qed{{\hskip0pt\unskip\unskip\nobreak\hfil\penalty50
          \hskip1em\hbox{}\nobreak\hfil
           {$\square$}
          \parfillskip=0pt\finalhyphendemerits=0
          \par}\medskip}

\newenvironment{Proof}
               {\noindent{\bf Proof.}\ }
               {\qed}

\newenvironment{Proofof}[1]
               {\noindent{\bf Proof of #1.}\ }
               {\qed}

\newcommand{\BA}{{\mathbb{A}}}

\newcommand{\BC}{{\mathbb{C}}}

\newcommand{\BF}{{\mathbb{F}}}
\newcommand{\BG}{{\mathbb{G}}}

\newcommand{\BP}{{\mathbb{P}}}

\newcommand{\Fm}{{\mathfrak{m}}}

\newcommand{\Fp}{{\mathfrak{p}}}

\newcommand{\CF}{{\cal F}}

\newcommand{\CL}{{\cal L}}
\newcommand{\CM}{{\cal M}}

\newcommand{\CO}{{\cal O}}

\newcommand{\CR}{{\cal R}}

\newbox\mybox
\def\arrover#1{\mathrel{
       \setbox\mybox=\hbox spread 1.4em
              {\hfil$\scriptstyle#1$\hfil}
       \vbox{\offinterlineskip\copy\mybox
             \hbox to\wd\mybox{\rightarrowfill}}}}

\def\larrover#1{\mathrel{
       \setbox\mybox=\hbox spread 1.4em
              {\hfil$\scriptstyle#1\vphantom{g}$\hfil}
       \vbox{\offinterlineskip\copy\mybox
             \hbox to\wd\mybox{\leftarrowfill}}}}

\def\ontoover#1{\mathrel{
       \setbox\mybox=\hbox spread 1.4em
              {\hfil$\scriptstyle#1\vphantom{g}$\hfil}
       \vbox{\offinterlineskip\copy\mybox
             \hbox to\wd\mybox{\rightarrowfill\hskip-2.8mm
                               $\rightarrow$}}}}
\def\leftontoover#1{\mathrel{
       \setbox\mybox=\hbox spread 1.4em
              {\hfil$\scriptstyle#1\vphantom{g}$\hfil}
       \vbox{\offinterlineskip\copy\mybox
             \hbox to\wd\mybox{$\leftarrow$\hskip-2.8mm
                               \leftarrowfill}}}}
\let\longto\longrightarrow
\let\into\hookrightarrow
\let\onto\twoheadrightarrow
\def\longonto{\ontoover{\ }}

\def\Cinf{{\BC}_\infty}
\def\Finf{F_\infty}

\def\Ahat{\hat{A}}

\begin{document}

\title{Compactification of Drinfeld Modular Varieties\\ and 
Drinfeld Modular Forms of Arbitrary Rank}

\author{Richard Pink$^{1,2}$}

\footnotetext[1]{Dept. of Mathematics, ETH Z\"urich, 8092 Z\"urich, Switzerland, {\tt pink@math.ethz.ch}}
\footnotetext[2]{Supported through the program ``Research in Pairs'' by Mathematisches Forschungsinstitut Oberwolfach in 2010.}

\date{March 16, 2012}
\maketitle

\begin{abstract}
We give an abstract characterization of the Satake compactification of a general Drinfeld modular variety. We prove that it exists and is unique up to unique isomorphism, though we do not give an explicit stratification by Drinfeld modular varieties of smaller rank which is also expected. We construct a natural ample invertible sheaf on it, such that the global sections of its $k$-th power form the space of (algebraic) Drinfeld modular forms of weight~$k$. We show how the Satake compactification and modular forms behave under all natural morphisms between Drinfeld modular varieties; in particular we define Hecke operators. We give explicit results in some special cases.
\end{abstract}

% \newpage
% \tableofcontents

%%%%%%%%%%%%%%%%%%%%%%%%%%%%%%%%%%%%%%%%%%%%%%%%%%%%%%%%%%%%%%%%%%%%%%%%%%%%%%%%%%%%%%%%%%%%%

% \newpage
\addtocounter{section}{-1}
\section{Introduction}
\label{Intro}

The theory of Drinfeld modular curves and Drinfeld modular forms of rank two is well-developed with a range of general and explicit results; see for instance Gekeler \cite{GekelerDMC} and Goss \cite{GossAUHP}, \cite{GossES}, \cite{GossMFFRT}. The aim of this article is to lay some groundwork for an algebro-geometric theory of modular forms on Drinfeld moduli spaces of arbitrary rank. It concentrates on the algebraic aspects of this topic, while a joint article planned with Breuer will deal with the analytic aspects and the translation between the two.

\medskip
{}From the point of view of algebraic geometry, a modular form of weight $k$ on any modular variety $M$ can be viewed as a section of the $k^{\rm th}$ power of a certain natural ample invertible sheaf on~$M$. When $M$ is a fine moduli space of Drinfeld modules, this invertible sheaf is the dual of the 
% sheaf of sections of the line bundle underlying
relative Lie algebra of the universal family of Drinfeld modules over~$M$. However, the fact that Drinfeld modular varieties are affine means that there is no K\"ocher principle, i.e., the definition of modular forms requires a condition at infinity. Thus an algebro-geometric definition of modular forms requires an extension of the invertible sheaf to a compactification $\OM$ of~$M$,
% (or at least a partial compactification), 
so that the space of sections over $\OM$ becomes finite dimensional. The natural candidate for $\OM$ is the analogue of the Satake compactification of Siegel moduli space, and the extension of the invertible sheaf to $\OM$ should arise naturally from the reduction of the universal family at the boundary $\OM\setminus M$.

\medskip
This sets the program for the present article: Describe $\OM$ and the behavior of the universal family at the boundary well enough to define the correct extension of the invertible sheaf and thus the space of algebraic modular forms of any integral weight $k$ on~$M$. Furthermore, do this functorially under all natural morphisms between Drinfeld modular varieties and the associated maps between modular forms, in particular under Hecke operators.

\medskip
It is expected, and has been proved by Kapranov \cite{Kapranov} in certain cases, that the Satake compactification of a Drinfeld modular variety of rank $r$ possesses a natural stratification by Drinfeld modular varieties of all ranks $\le r$ and can be constructed explicitly by piecing together quotients of Drinfeld period domains by rigid analytic means. Fortunately, the present program does not require such strong results in general and succeeds in a relatively pedestrian way. 

\medskip
Namely, we define a \emph{generalized Drinfeld module of rank $\le r$} over a scheme $S$ in the same way as a Drinfeld module of rank~$=r$, except that the fibers are required only to be Drinfeld modules of some, possibly varying, positive rank $\le r$. We call a generalized Drinfeld module \emph{weakly separating} if any isomorphism class of Drinfeld modules occurs in at most finitely many distinct fibers. For example, when $M$ is a sufficiently fine modular variety and thus possesses a universal family of Drinfeld modules, this universal family, stripped of its level structure, is weakly separating, because any fixed Drinfeld module possesses only finitely many level structures of a given type. We then characterize the Satake compactification $\OM$ axiomatically as any normal integral proper algebraic variety containing $M$ as an open dense subvariety, such that the universal family over $M$ extends to a weakly separating generalized Drinfeld module over~$\OM$. One of our main results, Theorem \ref{SC:SatExUni}, states that such a compactification always exists, is unique and projective, and that the extended generalized Drinfeld module is also unique. 

Then we define $\CL$ as the dual of the 
% sheaf of sections of the line bundle underlying
relative Lie algebra of 
the extended family over~$\OM$, and the space of \emph{algebraic modular forms of weight $k$} on $M$ as the space of global sections $H^0(\OM,\CL^k)$. We are convinced that this space corresponds to the space of analytically defined holomorphic modular forms of weight $k$ from Goss \cite[Defs.$\,$1.14, 1.54]{GossES}, but leave it to another article to carry out this identification.

We show that all the objects constructed behave in the expected way under the natural morphisms between Drinfeld modular varieties. In particular we describe the natural action of Hecke operators on Drinfeld modular forms. In the cases of small level where a universal family does not exist, we define the Satake compactification by taking quotients and define the space of modular forms by taking invariants under suitable automorphism groups.

In the last two sections we describe at large the special case $A=\BF_q[t]$ with level structure $(t)$, as well as certain quotients thereof. These results rely on a detailed study of the geometry of a compactification of a certain `finite Drinfeld period domain' that was carried out with Schieder in~\cite{PinkSchieder}. In particular we prove an observation of Breuer that the coefficients of the universal Drinfeld $\BF_q[t]$-module of rank $r$ form algebraically independent generators of the ring of modular forms of rank $r$ and level $1$ with respect to $\BF_q[t]$.

\medskip
The results of the present article can also be used to describe how the degree of a subvariety of $\OM$ behaves under Hecke operators. This is being applied in Hubschmid's Ph.~D.\ thesis \cite{Hubschmid} and may lead to simplifications in Breuer \cite{BreuerAO}.

\medskip
This article grew out of a joint project with Florian Breuer. It has profited from this collaboration in more ways than can be mentioned and would not exist without him. It is my pleasure to express my sincere gratitude to him. We are also grateful to the referee for pointing out a subtlety in the definition of Drinfeld modules over schemes that we address in Section \ref{GDM}.

%%%%%%%%%%%%%%%%%%%%%%%%%%%%%%%%%%%%%%%%%%%%%%%%%%%%%%%%%%%%%%%%%%%%%%%%%%%%%%%%%%%%%%%%%%%%%

\section{Drinfeld modular varieties}
\label{DMV}

Let $\BF_p$ denote the prime field of characteristic $p>0$. Let $F$ be a global function field of characteristic~$p$, that is, a finitely generated field of transcendence degree $1$ over~$\BF_p$.
Let $\infty$ be a fixed place of~$F$ with completion~$\Finf$, and let $\Cinf$
% = \hat{\bar{F}}_{\infty}$ 
be the completion of an algebraic closure of~$\Finf$.
Let $A$ denote the ring of elements in $F$ that are regular away from~$\infty$. 
Let $\hat{A} \cong \prod_{\Fp}A_\Fp$ be its profinite completion, and $\BA^f_F = \hat{A}\otimes_A F$ the ring of finite ad\`eles of $F$.

\medskip
Let $r$ be a positive integer, let $N$ be a non-zero proper ideal of~$A$, and let $S$ be a scheme over~$F$. A \emph{Drinfeld $A$-module of rank $r$} over $S$ is a pair $(E,\phi)$ consisting of a line bundle $E$ over $S$ and a ring homomorphism $\phi\!: A \to \End(E)$, $a\mapsto \phi_a$ satisfying the usual conditions (see Section \ref{GDM}). A (full) \emph{level $N$ structure} on it is an $A$-linear isomorphism of group schemes over~$S$
$$\lambda\!:\ \underline{(N^{-1}/A)^r} \stackrel\sim\longto \phi[N] := \bigcap_{a\in N} \Ker(\phi_a),$$
where $\underline{(N^{-1}/A)^r}$ denotes the constant group scheme over $S$ with fibers $(N^{-1}/A)^r$. Let
$$K(N) := \ker\big(\GL_r(\hat{A})\to\GL_r(A/N)\big)$$
denote the principal congruence subgroup of level~$N$. Let $M^r_{A,K(N)}$ denote the fine moduli space over $F$ of Drinfeld $A$-modules of rank $r$ with a level structure of level~$N$. This is an irreducible smooth affine algebraic variety of dimension $r-1$ of finite type over~$F$. 

(In fact, Drinfeld  defined a more general moduli functor over $\Spec A$ and proved that it is representable by an irreducible smooth affine scheme of relative dimension $r-1$ over $\Spec A$, provided that $N$ be contained in at least two distinct maximal ideals of~$A$: see \hbox{\cite[\S5]{Drinfeld1}}, \cite[Thm.$\,$1.8]{GekelerOhio}. Over $F$ one maximal ideal is enough, which is all we consider in this paper.)
% Note that this is only the generic fiber of the corresponding moduli scheme over $\Spec A$. The existence of this full moduli scheme --- see \cite[\S5]{Drinfeld1}, \cite[Thm.$\,$1.8]{GekelerOhio} --- requires that $N$ be contained in at least two distinct maximal ideals of~$A$, but over $F$ one maximal ideal is enough.) 
% It represents the functor $\CF_N : \Sch_F \longto \Set$ which to any scheme $S$ over $F$ associates the set of isomorphism classes of Drinfeld $A$-modules of rank $r$ with full level-$N$ structures.

% The action of $\GL_r(A/N)$ on $(N^{-1}/A)^r$ induces a natural action of $\GL_r(A/N)$, and hence of $\GL_r(\hat{A})$, on $M^r_{A,K(N)}$.

Consider another non-zero ideal $N'\subset N$ of~$A$. Then by restriction to the subgroup scheme $\underline{(N^{-1}/A)^r}$ of $\underline{(N^{\prime-1}/A)^r}$, any level $N'$ structure on a Drinfeld $A$-module of rank $r$ induces a level $N$ structure. This corresponds to a natural morphism of the moduli schemes
\addtocounter{Thm}{1}
\begin{equation}
\label{DMV:R1}
J_1 \!:\ M^r_{A,K(N')}\longto M^r_{A,K(N)}.
\end{equation}
Letting $K(N) \subset \GL_r(\hat{A})$ act on $M^r_{A,K(N')}$ through its action on level $N'$ structures, this induces an isomorphism
\addtocounter{Thm}{1}
\begin{equation}
\label{DMV:Quot}
M^r_{A,K(N')}/K(N) \stackrel{\sim}{\longto} M^r_{A,K(N)}.
\end{equation}
% (compare Section~\ref{MI}). 

For an arbitrary open compact subgroup $K\subset\GL_r(\Ahat)$ take any ideal $N$ as above such that $K(N) \subset K$. Then the action of $K$ on level $N$ structures induces an action on $M^r_{A,K(N)}$, and the isomorphy (\ref{DMV:Quot}) implies that the quotient
\addtocounter{Thm}{1}
\begin{equation}
\label{DMV:DMVDef}
M^r_{A,K} := M^r_{A,K(N)}/K
\end{equation}
is, up to a natural isomorphism, independent of the choice of~$N$. This is the \emph{Drinfeld modular variety of level~$K$}.

\begin{Def}
\label{DMV:FineDef}
A subgroup $K \subset \GL_r(\hat A)$ is called \emph{fine} if, for some maximal ideal $\Fp\subset A$,  the image of $K$ in $\GL_r(A/\Fp)$ is unipotent.
\end{Def}

Let $(E,\phi,\lambda)$ denote the universal family on $M^r_{A,K(N)}$. For $N'\subset N$ the universal family on $M^r_{A,K(N')}$ is the pullback of $(E,\phi,\lambda)$ under $J_1$ extended to a level $N'$ structure.

\begin{Prop}
\label{DMV:FineProp}
The action on $M^r_{A,K(N)}$ of any fine subgroup $K \subset \GL_r(\hat A)$ containing $K(N)$ factors through a free action of $K/K(N)$, and the family $(E,\phi)$ on $M^r_{A,K(N)}$ descends to a family of Drinfeld $A$-modules on $M^r_{A,K}$, which is independent of~$N$.
\end{Prop}

% \label{change_of_level}
% the morphism $J_1\!: M^r_{A,K(N)}\to M^r_{A,K}$ is an \`etale Galois covering with group $K/K(N)$, 

\begin{Proof}
A fixed point under $k\in K$ corresponds to a fiber $(E_x,\phi_x,\lambda_x)$ that is isomorphic to $(E_x,\phi_x,\lambda_x\circ k)$. This requires an automorphism $\xi$ of $E_x$ such that $\xi\circ\lambda_x=\lambda_x\circ k$. But Definition \ref{DMV:FineDef} implies that $N\subset\Fp$ and that the level structure $\lambda_x$ includes a non-zero $\Fp$-torsion point which is fixed by~$K$. Thus $\xi$ fixes that point and is therefore the identity. Now $\lambda_x=\lambda_x\circ k$ implies that $k\in K(N)$, proving the first assertion. The remaining assertions follow from the first and the remarks preceding the proposition.
\end{Proof}

\medskip
For simplicity we call the family of Drinfeld $A$-modules from Proposition \ref{DMV:FineProp} the \emph{universal family} on $M^r_{A,K}$. In fact, endowed with a certain additional structure it becomes the universal family making $M^r_{A,K}$ a fine moduli scheme, but we do not need this here. A consequence of (\ref{DMV:DMVDef}) and Proposition \ref{DMV:FineProp} is:

\begin{Prop}
\label{DMV:ModQuot}
For any open compact subgroup $K\subset\GL_r(\hat A)$ and any open normal subgroup $K'\triangleleft K$ we have $M^r_{A,K} \cong M^r_{A,K'}/K$. If moreover $K$ is fine, then $M^r_{A,K'} \onto M^r_{A,K}$ is an \'etale Galois covering with Galois group $K/K'$.
\end{Prop}

Drinfeld moduli spaces possess the following rigid analytic description. Let $\Omega^r$ denote the Drinfeld period domain obtained by removing all $F$-rational hyperplanes from the rigid analytic space $\BP^{r-1}(\Cinf)$. We can view it as the space of $\Cinf^\times$-equivalence classes of embeddings $\omega\!: F^r\into\Cinf$ whose $F_\infty$-linear extension $F_\infty^r\into\Cinf$ is still injective. Thus it carries a natural left action of 
$\GL_r(F)$ by setting $\gamma\omega := \omega\circ\gamma^{-1}$, and there is a natural isomorphism of rigid analytic spaces
\addtocounter{Thm}{1}
\begin{equation}
\label{DMV:DQ}
\GL_r(F)\!\bigm\backslash\! \bigl( \Omega^r\times\GL_r(\BA^f_F)/K \bigr)
\ \stackrel{\sim}{\longto}\ M^r_{A,K}(\Cinf).
\end{equation}
In the case $K=K(N)$ this isomorphism sends the equivalence class of a pair $(\omega,g) \in \Omega^r\times\GL_r(\BA^f_F)$ to the isomorphism class of the Drinfeld module with the lattice $\Lambda := \omega(F^r\cap g\hat A^r) \subset \Cinf$ and the level structure which makes the following diagram commute:
$$\xymatrix@C+20pt{
(N^{-1}/A)^r \ar[d]^{\wr}_{\rm id} \ar[rr]^-{\sim} && N^{-1}\Lambda/\Lambda \\
N^{-1}\hat A^r/\hat A^r \ar[r]^-{\sim}_-{g} &
N^{-1}g\hat A^r/g\hat A^r &
N^{-1}(F^r\cap g\hat A^r)/(F^r\cap g\hat A^r) \ar[l]^-{\rm id}_-{\sim} \ar[u]_\omega^\wr \\}$$
% Moreover, the left hand side of (\ref{DMV:DQ}) can be identified with the disjoint union of quotients $\Gamma_x \backslash \Omega^r$, where $x$ runs through a set of representatives for $\GL_r(F)\backslash \GL_r(\BA_F^f)/K$ and $\Gamma_x = \GL_r(F)\cap x^{-1}Kx$. If $K$ is fine, then all $\Gamma_x$ act without fixed points on $\Omega^r$, so that the quotient map $\Omega^r \onto \Gamma_x \backslash \Omega^r$ is a local isomorphism everywhere.
 
The left hand side of (\ref{DMV:DQ}) makes sense for all open compact subgroups $K\subset\GL_r(\BA_F^f)$, not necessarily contained in $\GL_r(\hat A)$. Though one can extend the definition of $M^r_{A,K}$ accordingly, we refrain from doing so, because by conjugating $K$ back into $\GL_r(\hat A)$ one can identify these new spaces with previous ones that possess a more natural modular interpretation.

%%%%%%%%%%%%%%%%%%%%%%%%%%%%%%%%%%%%%%%%%%%%%%%%%%%%%%%%%%%%%%%%%%%%%%%%%%%%%%%%%%%%%%%%%%%%%

\section{Morphisms and isogenies}
\label{MI}

Consider an element $g\in\GL_r(\BA_F^f)$ with coefficients in~$\hat A$, so that $\hat A^r \subset g^{-1}\hat A^r$. Consider non-zero ideals $N'\subset N\subsetneqq\nobreak A$ such that $g^{-1}\hat A^r \subset N^{\prime-1}N\hat A^r$. Then we have a short exact sequence
$$\xymatrix{0 \ar[r] & g^{-1}\hat A^r / \hat A^r \ar[r] &
N^{-1}g^{-1}\hat A^r / \hat A^r \ar[r]^-g &
N^{-1}\hat A^r / \hat A^r \ar[r] & 0,\\}$$
where the middle term is contained in $N^{\prime-1}\hat A^r / \hat A^r \cong (N^{\prime-1}/A)^r$ and the right hand term is isomorphic to $(N^{-1}/A)^r$.
Thus for any Drinfeld $A$-module with level $N'$ structure $(E',\phi',\lambda')$ we can form the quotient by the constant torsion subgroup $(E',\phi')/\lambda'(g^{-1}\hat A^r/\hat A^r)$ and endow it with the level $N$ structure $\lambda$ defined by the formula $\lambda(gm) = \lambda'(m) \mod \lambda'(g^{-1}\hat A^r/\hat A^r)$ for all $m\in N^{-1}g^{-1}\hat A^r/\hat A^r$. Using the modular interpretation this defines a morphism of Drinfeld modular varieties 
\addtocounter{Thm}{1}
\begin{equation}
\label{MI:R1a}
J_g\!:\ M^r_{A,K(N')} \to M^r_{A,K(N)}.
\end{equation}
% If $(E',\phi',\lambda')$ denotes the universal family on $M^r_{A,K(N')}$ and $(E,\phi,\lambda)$ the universal family on $M^r_{A,K(N)}$, the modular characterization of $J_g$ means that $J_g^*(E,\phi,\lambda)$ is obtained from
If $(E,\phi,\lambda)$ denotes the universal family on $M^r_{A,K(N)}$, the modular characterization of $J_g$ means that $J_g^*(E,\phi,\lambda)$ is obtained from the universal family $(E',\phi',\lambda')$ on $M^r_{A,K(N')}$ by the above prescription. In particular, we obtain an isogeny
\addtocounter{Thm}{1}
\begin{equation}
\label{MI:R3}
\xi_g\!:\ (E',\phi') \longto J_g^*(E,\phi)
% (E,\phi|A',\lambda\circ e)\ \cong\ I_e^*(E',\phi',\lambda').
\end{equation}
with kernel $\lambda'(g^{-1}\hat A^r/\hat A^r)$. 

In the case $g=1$ the morphism $J_1$ from (\ref{MI:R1a}) is the same as that from (\ref{DMV:R1}), and the isogeny (\ref{MI:R3}) is an isomorphism. The following proposition gives another example:

\begin{Prop}
\label{MI:Scalar}
For any non-zero scalar $a\in A$ and any non-zero ideals $N\subsetneqq\nobreak A$ and $N'\subset aN$ we have $J_a=J_1$ and $\xi_a = \xi_1 \circ\phi'_a = \phi_a\circ\xi_1$.
\end{Prop}

\begin{Proof}
By construction $J_a$ sends $(E',\phi',\lambda')$ to the Drinfeld $A$-module $(E',\phi')/\lambda'(a^{-1}\hat A^r/\hat A^r)$ with the level $N$ structure $\lambda$ that is characterized by $\lambda(am) = \lambda'(m) \mod \lambda'(a^{-1}\hat A^r/\hat A^r)$. But $\lambda'(a^{-1}\hat A^r/\hat A^r)$ is the kernel of the isogeny $\phi'_a\!: (E',\phi')\to (E',\phi')$; hence $\phi'_a$ induces an isomorphism $(E',\phi')/\lambda'(a^{-1}\hat A^r/\hat A^r) \stackrel{\sim}{\longto} (E',\phi')$. Under this isomorphism, the induced level $N$ structure $\lambda$ corresponds to the homomorphism $am\mapsto \phi'_a(\lambda'(m)) = \lambda'(am)$, which is simply the restriction $\lambda'|(N^{-1}/A)^r$. The resulting data is thus isomorphic to that obtained by~$J_1$, and everything follows.
\end{Proof}

\medskip
Now consider a second element $g'\in\GL_r(\BA_F^f)$ with coefficients in~$\hat A$ and a third non-zero ideal $N''\subset N'$ such that $g^{\prime-1}\hat A^r \subset N^{\prime\prime-1}N'\hat A^r$, and let $(E'',\phi'',\lambda'')$ be the universal family on $M^r_{A,K(N'')}$. Then a direct calculation shows that
\addtocounter{Thm}{1}
\begin{equation}
\label{MI:R2xa}
J_g\circ J_{g'} = J_{gg'}
\end{equation}
and that $\xi_{gg'}$ is the composite of the isogenies
\addtocounter{Thm}{1}
\begin{equation}
\label{MI:R4}
(E'',\phi'') \stackrel{\xi_{g'}}{\longto} J_{g'}^*(E',\phi')
\arrover{J_{g'}^*\xi_g} J_{g'}^*J_g^*(E,\phi) = J_{gg'}^*(E,\phi).
\end{equation}

%%%%%%%%%%%%%%%%%%%%%%%

For any element $g\in\GL_r(\hat A)$ we can take $N'=N$ and obtain an automorphism $J_g$ of $M^r_{A,K(N)}$ and a covering automorphism $\xi_g$ of $(E,\phi)$. The relation (\ref{MI:R2xa}) shows that this defines a left action of $\GL_r(\hat A)$. This is precisely the action used in Section \ref{DMV}
% Proposition~\ref{DMV:FineProp} 
and helps to extend the above constructions to more general open compact subgroups, as follows:

\begin{Prop}
\label{MI:FineProp1}
Consider an element $g\in\GL_r(\BA_F^f)$ and two open compact subgroups $K$, $K'\subset\GL_r(\hat A)$ such that $gK'g^{-1}\subset K$. Then there is a natural morphism 
$$J_g\!:\ M^r_{A,K'} \to M^r_{A,K},$$
which coincides with that in (\ref{MI:R1a}) if $K'=K(N')$ and $K=K(N)$ under the assumptions there. 
If $g'\in\GL_r(\BA_F^f)$ is a second element and $K''\subset\GL_r(\hat A)$ a third open compact subgroup such that $g'K''g^{\prime-1}\subset K'$, then 
$$J_g\circ J_{g'} = J_{gg'}.$$
Furthermore, the morphism $J_g$ is finite.
\end{Prop}

\begin{Proof}
Assume first that $g$ has coefficients in~$\hat A$, so that $\hat A^r \subset g^{-1}\hat A^r$. Choose a non-zero proper ideal $N\subsetneqq\nobreak A$ such that $K(N)\subset K$. Thereafter choose a non-zero ideal $N'\subset N$ such that $K(N')\subset K'$ and $g^{-1}\hat A^r \subset N^{\prime-1}N\hat A^r$. Then all the assumptions on $(g,N',N)$ are satisfied for the morphism $J_g$ from (\ref{MI:R1a}). Consider the composite morphism 
$$M^r_{A,K(N')} \stackrel{J_g}{\longto} M^r_{A,K(N)} \longonto M^r_{A,K(N)}/K\ =\ M^r_{A,K}.$$
Using (\ref{MI:R2xa}) one easily shows that this morphism is invariant under the action of $K'$ on $M^r_{A,K(N')}$; hence it factors through a unique morphism $J_g\!:\ M^r_{A,K'} \to M^r_{A,K}$. Direct calculations show that this morphism is independent of the choice of $N$ and~$N'$.
Moreover, if $g'$ and hence $gg'$ also has coefficients in~$\hat A$, the relation (\ref{MI:R2xa}) implies that the new morphisms thus obtained satisfy the relation $J_g\circ J_{g'} = J_{gg'}$. In particular, for any non-zero scalar $a\in A$ Proposition \ref{MI:Scalar} implies that these morphisms satisfy $J_g = J_1\circ J_g = J_a\circ J_g = J_{ag}$.

For arbitrary $g$ consider any non-zero scalar $a\in A$ such that $ag$ has coefficients in~$\hat A$. Then the relation just proved implies that the morphism $J_{ag}\colon M^r_{A,K'} \to M^r_{A,K}$ is independent of~$a$. We can therefore generally define $J_g:= J_{ag}$ for any non-zero scalar $a\in A$ such that $ag$ has coefficients in~$\hat A$. By a short calculation these morphisms inherit the relation $J_g\circ J_{g'} = J_{gg'}$. This proves the first two assertions of the proposition.

To prove the last assertion, by the construction of $J_g$ we may assume that $g$ has coefficients in~$\hat A$. Fix a non-zero element $a\in A$ such that $h:= ag^{-1}$ has coefficients in~$\hat A$. Take any non-zero ideal $N_0\subsetneqq A$, and in the above choice of $N$ and $N'$ assume in addition that $h^{-1}\hat A^r \subset N^{-1}N_0\hat A^r$. Then we have morphisms
$$M^r_{A,K(N')} \stackrel{J_g}{\longto} M^r_{A,K(N)} \stackrel{J_h}{\longto} M^r_{A,K(N_0)}$$
whose composite is $J_h\circ J_g=J_{hg}=J_a$. By Proposition \ref{MI:Scalar} this is equal to $J_1$ and hence finite by \cite[\S5]{Drinfeld1}, \cite[Thm.$\,$1.8]{GekelerOhio}. As all varieties are separated, using \cite[Ch.$\,$II Cor.$\,$4.8 (e)]{Hartshorne} it follows that $J_g\!: M^r_{A,K(N')} \to M^r_{A,K(N)}$ is finite. Since the morphism $J_g\!: M^r_{A,K'} \to M^r_{A,K}$ is obtained from this by taking quotients, it is also finite, as desired.
\end{Proof}

\begin{Prop}
\label{MI:FineProp2}
In Proposition \ref{MI:FineProp1} assume moreover that $g$ and $g'$ have coefficients in $\hat A$ and that $K$, $K'$, $K''$ are fine. Let $(E,\phi)$, $(E',\phi')$, $(E'',\phi'')$ denote the respective universal families on $M^r_{A,K}$, $M^r_{A,K'}$, $M^r_{A,K''}$. Then there is a natural isogeny
$$\xi_g\!:\ (E',\phi') \longto J_g^*(E,\phi),$$
which coincides with that in (\ref{MI:R3}) if $K'=K(N')$ and $K=K(N)$ under the assumptions there. 
Furthermore these isogenies satisfy 
$$(J_{g'}^*\xi_g)\circ \xi_{g'} = \xi_{gg'}.$$
\end{Prop}

\begin{Proof}
(Sketch) 
Let $N$ and $N'$ be as in the proof of Proposition \ref{MI:FineProp1}. Then (\ref{MI:R3}) yields an isogeny $J_1^*(E',\phi') \to J_1^*J_g^*(E,\phi)$ after pullback via $J_1\!: M^r_{A,K(N')} \onto M^r_{A,K'}$. But the cocycle relation (\ref{MI:R4}) implies that this isogeny over $M^r_{A,K(N')}$ is equivariant under the action of $K'/K(N')$. It therefore descends to an isogeny over $M^r_{A,K'}$. Using the equivariance one also shows that the isogeny thus constructed is independent of $N$ and~$N'$. This proves the first assertion. The second assertion follows from (\ref{MI:R4}) by a similar calculation.
\end{Proof}

%%%%%%%%%%%%%%%%%%%%%%%

\medskip
We leave it to the careful reader to verify that under the isomorphism (\ref{DMV:DQ}) the morphism $J_g$ from Proposition \ref{MI:FineProp1} corresponds to the morphism of rigid analytic spaces induced by right translation
\addtocounter{Thm}{1}
\begin{eqnarray}
\label{MI:R1}
% 	M^{r'}_{A',K'} & \stackrel{J_g}{\longto} & M^r_{A,K}  \\
% J_g\!:\ \ 
	\GL_r(F)\!\bigm\backslash\! \bigl( \Omega^r\times\GL_r(\BA^f_F)/K' \bigr) &\longto &
	\GL_r(F)\!\bigm\backslash\! \bigl( \Omega^r\times\GL_r(\BA^f_F)/K \bigr), \\
\nonumber	\;[(\omega,h)] & \longmapsto & [(\omega,hg^{-1})].
\end{eqnarray}

%%%%%%%%%%%%%%%%%%%

\medskip
To describe morphisms between Drinfeld modular varieties of different types, let $F'$ be a finite extension of $F$ which possesses a unique place $\infty'$ above the place $\infty$ of~$F$. Then the ring $A'$ of elements of $F'$ that are regular away from $\infty'$ is the integral closure of $A$ in~$F'$. Assume that $r = r'\cdot [F'/F]$ for a positive integer~$r'$, and choose an $\hat A$-linear isomorphism $b\!: \hat A^r \stackrel{\sim}{\to} \hat A^{\prime\kern1pt r'}$. Take a non-zero proper ideal $N\subsetneqq A$ and set $N':=NA'$. Then $b$ induces an $A$-linear isomorphism $(N^{-1}/A)^r \stackrel\sim\to (N^{\prime-1}/A')^{r'}$ that we again denote by~$b$. We can thus associate to any Drinfeld $A'$-module with level $N'$ structure $(E',\phi',\lambda')$ the Drinfeld $A$-module with level $N$ structure $(E',\phi'|A,\lambda'\circ b)$. Using the modular interpretation this defines a morphism of Drinfeld modular varieties
\addtocounter{Thm}{1}
\begin{equation}
\label{MI:Theta1}
I_b\!:\ M^{r'}_{A',K(N')} \longto M^r_{A,K(N)}.
\end{equation}
By \cite[Lem.$\,$3.1, Prop.$\,$3.2]{BreuerAO} this morphism is injective and finite. If $(E',\phi',\lambda')$ denotes the universal family on $M^{r'}_{A',K(N')}$ and $(E,\phi,\lambda)$ the universal family on $M^r_{A,K(N)}$, this morphism is characterized uniquely by an isomorphism
\addtocounter{Thm}{1}
\begin{equation}
\label{MI:Theta2}
(E',\phi'|A,\lambda'\circ b)\ \cong\ I_b^*(E,\phi,\lambda).
\end{equation}

More generally, for arbitrary open compact subgroups $K'\subset\GL_{r'}(\hat A')$ and $K\subset\GL_r(\hat A)$ choose $N$ and $N':=NA'$ such that $K(N') \subset K'$ and $K(N) \subset K$. Then the composite morphism $J_1\circ I_b\!: M^{r'}_{A',K(N')} \to M^r_{A,K(N)} \onto M^r_{A,K}$ is $K'$-invariant provided that $K'$ is mapped into $K$ under the embedding 
$$\GL_{r'}(\hat A') \into \GL_r(\hat A),\ \ k'\mapsto b^{-1}k'b.$$ 
In this case it factors through a finite morphism of Drinfeld modular varieties
\addtocounter{Thm}{1}
\begin{equation}
\label{MI:Theta3}
I_b\!:\ M^{r'}_{A',K'} \longto M^r_{A,K}.
\end{equation}

For a rigid analytic description of~$I_b$ consider in addition an $F$-linear isomorphism $\beta\!: 
F^r \stackrel{\sim}{\to}\nobreak F^{\prime\kern1pt r'}$.  Then both $b$ and $\beta$ induce isomorphisms $\BA_F^r \stackrel{\sim}{\to} \BA_{F'}^{r'}$ that we again denote by $b$ and~$\beta$. (They can be made to coincide if $A'$ is a free $A$-module, but not in general.) A direct calculation shows that the morphism (\ref{MI:Theta3}) corresponds to the map
\addtocounter{Thm}{1}
\begin{eqnarray}
\label{MI:Theta4}
\GL_{r'}(F')\!\bigm\backslash\! \bigl( \Omega^{r'}\times\GL_{r'}(\BA^f_{F'})/K' \bigr)
& \longto &
\GL_r(F)\!\bigm\backslash\! \bigl( \Omega^r\times\GL_r(\BA^f_F)/K \bigr), \\
\nonumber {}[(\omega',g')] & \longmapsto & \hfil[(\omega'\circ\beta,\beta^{-1}g'b)].
\end{eqnarray}
Note that the equivalence class on the right hand side is in any case independent of~$\beta$; the choice of $\beta$ is needed only to write down a representative for it.

%%%%%%%%%%%%%%%%%%%%%%%%%%%%%%%%%%%%%%%%%%%%%%%%%%%%%%%%%%%%%%%%%%%%%%%%%%%%%%%%%%%%%%%%%%%%%

\section{Generalized Drinfeld modules}
\label{GDM}

The definition of Drinfeld modules over a scheme involves a subtlety over which one can easily stumble, the present author included. The choice in dealing with that subtlety is not important when studying Drinfeld modules of constant rank, but it has a non-trivial effect for degenerating Drinfeld modules. We therefore take some time to discuss the notion in detail.

\medskip
%Let $S$ be a scheme over~$\BF_p$. By definition a \emph{line bundle over $S$} is a group scheme $E$ over $S$ together with a morphism $\BG_a\times E\to E$ which is Zariski locally over $S$ isomorphic to the additive group scheme $\BG_{a,S}$ together with the multiplication morphism $\BG_a\times \BG_{a,S}\to \BG_{a,S}$, $(x,y)\mapsto xy$.
By definition the \emph{trivial line bundle over} a scheme $S$ is the additive group scheme $\BG_{a,S}$ over $S$ together with the morphism $\BG_a\times \BG_{a,S}\to \BG_{a,S}$, $(x,y)\mapsto xy$. An arbitrary \emph{line bundle over $S$} is a group scheme $E$ over $S$ together with a \emph{scalar multiplication} morphism $\BG_a\times E\to E$ which, as a pair, is Zariski locally over $S$ isomorphic to the trivial line bundle. 
% Note that, although giving a line bundle is equivalent to giving the invertible sheaf of $\CO_S$-modules of its sections, we keep these notions separate.
The automorphism group of any line bundle $E$ over $S$ is $\BG_m(S)$, acting on $E$ through the given scalar multiplication.

By contrast one can consider just a group scheme over $S$ which is Zariski locally isomorphic to $\BG_{a,S}$. Any line bundle yields such a group scheme by forgetting the scalar multiplication, but it is not at all clear whether a group scheme which is Zariski locally isomorphic to $\BG_{a,S}$ comes from a line bundle or whether that line bundle is unique or at least unique up to isomorphism. In characteristic zero the answer to these questions is yes, but in positive characteristic the situation is different.

So let $S$ be a scheme over~$\BF_p$, for simplicity quasi-compact, and let $E$ be a line bundle over~$S$. Let $\End(E)$ denote the ring of endomorphisms of the (commutative) group scheme underlying~$E$. 
% In the case $E=\BG_{a,S}$ this is the non-commutative polynomial ring $\Gamma(S,\CO_S)[\tau]$ with generator $\tau\!: x\mapsto x^p$. 
As observed by Drinfeld \cite[\S5]{Drinfeld1}, any such endomorphism can be written uniquely as a finite sum $\sum_i b_i\tau^i$ for sections $b_i\in\Gamma(S,E^{1-p^i})$ and the Frobenius homomorphism $\tau: E\to E^{p}$, $x\mapsto x^p$. Such an endomorphism is an automorphism of group schemes if and only if $b_0\in\Gamma(S,E^0) = \Gamma(S,\CO_S)$ is invertible and $b_i$ is nilpotent for every $i>0$. It is an automorphism of the line bundle if and only if $b_0$ is invertible and $b_i=0$ for all $i>0$. 

Thus if $S$ is affine and not reduced, there exist automorphisms of the group scheme underlying $E$ which do not commute with the scalar multiplication. Twisting the scalar multiplication by such an automorphism then yields a different structure of line bundle on the same underlying group scheme~$E$. Consequently, if a group scheme over $S$ comes from a line bundle, that line bundle is in general not unique. It is therefore important to distinguish these two notions.

%(The reader is invited to prove that if a group scheme over $S$ comes from a line bundle~$E$, then $E$ is determined up to a natural isomorphism, because $E$ is naturally isomorphic to the line bundle associated to its relative Lie algebra, and this Lie algebra depends only on the underlying group scheme of~$E$. More precisely, if $E$ and $E'$ are line bundles on $S$ with the same underlying group scheme, there is a unique isomorphism of line bundles $f:E\to E'$ which induces the identity on the relative Lie algebra. However, $f$ is not necessarily the identity on the underlying group scheme.)

The following definitions are based on the notion of line bundles, not just group schemes locally isomorphic to $\BG_{a,S}$, for reasons explained below.

\medskip
Let $A$ be the ring from Section~\ref{DMV}. The \emph{degree} of a non-zero element $a\in A$ is the integer $\deg(a)\ge0$ defined by the equality $[A:(a)] = p^{\deg(a)}$. 

\medskip
First consider a line bundle $E$ on the spectrum of a field $K$ and a homomorphism $\phi\!: A \to \End(E)$ given by $a\mapsto \phi_a = \sum_i \phi_{a,i} \tau^i$ with $\phi_{a,i}$ in the one-dimensional $K$-vector space $\Gamma(\Spec K,E^{1-p^i})$. By Drinfeld \cite[Prop.$\,$2.1 \& Cor.]{Drinfeld1} or \cite[Prop.$\,$4.5.3]{GossBS}, there exists a unique integer $r\ge0$ such that $\phi_{a,i}=0$ whenever $i>r\deg(a)$ and $\phi_{a,r\deg(a)}\not=0$ whenever $r\deg(a)>0$. If this integer is $>0$, then $\phi$ is called a \emph{Drinfeld $A$-module of rank $r$ over~$K$}.

\medskip
Let $S$ be a scheme over $\Spec A$. 

\begin{Def}
\label{GDM:GenDrinDef}
A \emph{generalized Drinfeld $A$-module over~$S$} is a pair $(E,\phi)$ consisting of a line bundle $E$ over~$S$ and a ring homomorphism 
$$\textstyle \phi\!: A \to \End(E),\ a\mapsto \phi_a = \sum_i \phi_{a,i} \tau^i$$
with $\phi_{a,i}\in\Gamma(S,E^{1-p^i})$ satisfying the following two conditions:
\begin{itemize}
\item[(a)] The derivative $d\phi\!: a\mapsto \phi_{a,0}$ is the structure homomorphism $A\to\Gamma(S,\CO_S)$.
% The map $d\phi\!: A\to\Gamma(S,\CO_S)$ . . . 
\item[(b)] Over any point $s\in S$ the map $\phi$ defines a Drinfeld module of some rank $r_s\ge1$.
\end{itemize} 
A generalized Drinfeld $A$-module is \emph{of rank $\leq r$} if
\begin{itemize}
\item[(c)] For all $a\in A$ and $i>r\deg(a)$ we have $\phi_{a,i}=0$.
\end{itemize} 
An \emph{isomorphism} of generalized Drinfeld $A$-modules is an isomorphism of line bundles that is equivariant with respect to the action of $A$ on both sides.
\end{Def}

For any generalized Drinfeld module, the function $s\mapsto r_s$ is lower semicontinuous, because any coefficient which is non-zero at a point remains non-zero in a neighborhood. If the generalized Drinfeld module is of rank $\le r$, we have $\max\,\{r_s\,|\,s\in S\} \le r$.
% $r_s\le r$ for all~$s$. 
The converse is true if $S$ is reduced, because then a section of a line bundle on $S$ is zero if and only if its value at every point of $S$ is zero. In general, however, it is possible that a higher coefficient $\phi_{a,i}$ is nilpotent, and so a generalized Drinfeld module may not be of rank $\le \max\{r_s\mid s\in S\}$. In that case we can view it as a non-trivial infinitesimal deformation towards a Drinfeld module of higher rank, and our notion is geared precisely towards capturing this possibility. We hope that this will provide a better basis for some kind of modular interpretation of generalized Drinfeld modules. Note also that the notion `of rank $\le r$' is invariant under isomorphisms.

%%%%%%%%%%%%%%%%%%%%%%%
\medskip

\begin{Def}
\label{GDM:DrinDef}
A generalized Drinfeld $A$-module of rank $\leq r$ with $r_s=r$ everywhere is called a \emph{Drinfeld $A$-module of rank~$r$} over~$S$.
\end{Def}

% Thus a homomorphism $\phi\!: A \to \End(E)$, $a\mapsto \phi_a = \sum_i \phi_{a,i} \tau^i$ defines a Drinfeld module of rank $r\ge1$ if and only if condition \ref{GDM:GenDrinDef} (a) holds and $\phi_{a,i}=0$ whenever $i>r\deg(a)$ and $\phi_{a,r\deg(a)}$ is invertible whenever $r\deg(a)>0$.

\begin{Rem}
\label{GDM:DrinDefRem}
\rm This definition corresponds to that of a `standard' elliptic $A$-module from Drinfeld \cite[\S5B]{Drinfeld1}, which is suggested as a variant of the one officially used there,
and which was also adopted in \cite[Def.$\,$1.2]{Saidi}. The original definition in \cite[\S5B]{Drinfeld1} requires a generalized Drinfeld $A$-module with $r_s=r$ everywhere, without our condition \ref{GDM:GenDrinDef} (c), and an isomorphism of Drinfeld modules is defined there as any isomorphism of the underlying group schemes that is equivariant under~$A$. That notion, as it stands, does not lend itself to gluing over a Zariski open covering, because, although $E$ is required to be a line bundle, the possible gluing isomorphisms may not be isomorphisms of line bundles, and so the glued group scheme may not inherit a natural structure of line bundle. If one follows this approach, it would be more natural to 
% define a Drinfeld $A$-module of rank $r$ as a pair consisting only of a group scheme $E$ over $S$ which is locally isomorphic to $\BG_{a,S}$ together with a homomorphism $\phi\colon A\to\End(S)$ satisfying \ref{GDM:GenDrinDef} (a--b) and with $r_s=r$ everywhere.
% drop the line bundle altogether and replace it by 
replace the line bundle throughout by a group scheme over $S$ which is locally isomorphic to $\BG_{a,S}$, which would make the problem disappear.
The following fact, adapted from \cite[\S5B]{Drinfeld1}, shows that the resulting theory is equivalent to that using the above Definition~\ref{GDM:DrinDef}:
\end{Rem}

\begin{Prop}
\label{GDM:DrinDefProp}
Let $E$ be a group scheme over $S$ which is locally isomorphic to $\BG_{a,S}$, and let $\phi\!: A \to \End(E)$ be a homomorphism satisfying the conditions \ref{GDM:GenDrinDef} (a--b) with $r_s=r$ everywhere. Then $E$ possesses a unique structure of line bundle making $(E,\phi)$ into a Drinfeld $A$-module of rank $r$ according to Definition~\ref{GDM:DrinDef}.
\end{Prop}

\begin{Proof}
By uniqueness, it suffices to prove everything locally over~$S$. Thus we may assume that $S$ is affine 
% , say $S=\Spec R$, 
and that $E=\BG_{a,S}$ as a group scheme over~$S$. Choose any non-constant element $t\in A$. Then by \cite[Prop.$\,$5.2]{Drinfeld1}, there exists a unique automorphism $f$ of the group scheme $\BG_{a,S}$ which is the identity on the Lie algebra, such that 
$$f\phi_t f^{-1}\ =\ \sum_{i=0}^{r\deg(t)} u_i\cdot\tau^i$$
with $u_i\in R$ and $u_{r\deg(t)}\in R^\times$. For every non-constant element $a\in A$ it then follows from \cite[Prop.$\,$5.1]{Drinfeld1} and the fact that $f\phi_tf^{-1}$ and $f\phi_af^{-1}$ commute that 
$$f\phi_a f^{-1}\ =\ \sum_{i=0}^{r\deg(a)} v_i\cdot\tau^i$$
with $v_i\in R$ and $v_{r\deg(a)}\in R^\times$.
Thus the trivial line bundle $E_0 := \BG_{a,S}$ and the map $a\mapsto f\phi_a f^{-1}$ constitute a Drinfeld $A$-module of rank $r$ according to Definition~\ref{GDM:DrinDef}. Now, by transport of structure, the group scheme $E=\BG_{a,S}$ possesses a unique structure of line bundle such that $f$ induces an isomorphism of line bundles $E\stackrel{\sim}{\to}E_0$. With this structure, the pair $(E,\phi)$ is then a Drinfeld $A$-module of rank $r$ according to Definition~\ref{GDM:DrinDef}. This proves the existence part.

To prove the uniqueness consider any structure of line bundle on $E$ such that $(E,\phi)$ is a Drinfeld module of rank $r$ according to Definition~\ref{GDM:DrinDef}. 
After possibly localizing on $S$ we may assume that there exists an isomorphism of line bundles $g\colon E\stackrel{\sim}{\to}E_0$. Then $g$ acts on the common Lie algebra of the underlying group scheme $\BG_{a,S}$ through multiplication by a unit $u\in\Gamma(S,\CO_S^\times)$. After replacing $g$ by $u^{-1}g$ we may thus assume that $g$ induces the identity on the Lie algebra. Then $E_0$ and the map $a\mapsto g\phi_a g^{-1}$ constitute a Drinfeld $A$-module of rank $r$ according to Definition~\ref{GDM:DrinDef}. In particular we have 
$$g\phi_t g^{-1}\ =\ \sum_{i=0}^{r\deg(t)} w_i\cdot\tau^i$$
with $w_i\in R$ and $w_{r\deg(t)}\in R^\times$. But by the uniqueness of $f$ this implies that $g=f$. It follows that the structure of line bundle on $E$ is equal to that given by~$f$; hence it is unique, as desired.
\end{Proof}

%%%%%%%%%%%%%%%%%%%%%%%

\begin{Def}
\label{GDM:HomDef}
A \emph{homomorphism} $\xi\!: (E,\phi) \to (E',\phi')$ of generalized Drinfeld $A$-modules over $S$ is a homomorphism of the underlying group schemes $\xi\!:E\to E'$ satisfying $\xi\circ\phi_a = \phi'_a\circ\xi$ for all $a\in A$. A homomorphism which is non-zero in every fiber is called an \emph{isogeny}. 
\end{Def}

For example, any automorphism $f$ of the group scheme underlying $E$ determines an isogeny from $(E,\phi)$ to another generalized Drinfeld module $(E,f\phi f^{-1})$. By construction this isogeny has a two-sided inverse, though it may not necessarily be an isomorphism of generalized Drinfeld modules according to Definition \ref{GDM:GenDrinDef} if it is not also an automorphism of line bundles. This is an unfortunate drawback of the definition. At least the problem disappears in the following cases:

\begin{Prop}
\label{GDM:DrinIsomProp}
Let $\xi$ be a homomorphism of generalized Drinfeld modules over $S$ which possesses a two-sided inverse. If $S$ is reduced, or if both generalized Drinfeld modules are Drinfeld modules according to Definition \ref{GDM:DrinDef}, then $\xi$ is an isomorphism.
\end{Prop}

\begin{Proof}
Since $\xi$ has a two-sided inverse, it is an isomorphism of the group schemes underlying the given line bundles. If $S$ is reduced, any such isomorphism is already an isomorphism of line bundles. If the generalized Drinfeld modules are Drinfeld modules according to Definition \ref{GDM:DrinDef}, the same conclusion follows from the uniqueness in Proposition \ref{GDM:DrinDefProp}. In both cases the proposition follows.
\end{Proof}

%%%%%%%%%%%%%%%%%%%%%%%

\begin{Prop}
\label{GDM:HomExt}
Assume that $S$ is normal integral and $U\subset S$ is open dense. Let $(E,\phi)$ and $(E',\phi')$ be generalized Drinfeld $A$-modules over~$S$. Then any homomorphism $\xi_U\!: (E,\phi)|U \to (E',\phi')|U$ extends to a unique homomorphism $\xi\!: (E,\phi) \to (E',\phi')$.
\end{Prop}

\begin{Proof}
As the problem is local on $S$, we may assume that $S$ is the spectrum of a normal integral domain~$R$ and that $E=E'=\BG_{a,S}$. Let $K$ denote the quotient field of~$R$, and let $\Spec R'$ be a non-empty open affine in~$S$. Then $\phi$ and $\phi'$ are homomorphisms $A\to\End(\BG_{a,S}) = R[\tau]$, and $\xi_U$ is an element of $R'[\tau]$. We must show that $\xi_U$ actually lies in $R[\tau]$. Since $R$ is integrally closed, it is the intersection of all valuation rings containing it by \cite[Thm.$\,$10.4]{Matsumura}. Thus it suffices to prove the assertion after replacing $R$ by any valuation ring containing it. Let then $\Fm$ be the maximal ideal of~$R$.

Assume that $\xi_U \not\in R[\tau]$. Take an element $\lambda\in K\setminus R$ such that $\lambda^{-1}\xi_U \in R[\tau] \setminus\Fm[\tau]$; for instance a coefficient of $\xi_U$ of minimal valuation.
% $\lambda^{-1}\xi_U$ has all coefficients in~$R$ but not all coefficients in~$\Fm$. 
Fix any non-constant element $a\in A$. We claim that $\lambda^{-1}\phi'_a\lambda \in R[\tau]$.
Indeed, if that is not the case, take $\mu\in K\setminus R$ such that $\mu^{-1}\lambda^{-1}\phi'_a\lambda \in R[\tau] \setminus\Fm[\tau]$. Then the defining relation for $\xi_U$ implies that
$$(\mu^{-1}\lambda^{-1}\phi'_a\lambda)(\lambda^{-1}\xi_U)
\ =\ \mu^{-1}\lambda^{-1}\phi'_a\xi_U
\ =\ \mu^{-1}\lambda^{-1}\xi_U\phi_a
\ =\ \mu^{-1}(\lambda^{-1}\xi_U)\phi_a.$$
Here the left hand side is in $R[\tau] \setminus\Fm[\tau]$, because the ring $(R/\Fm)[\tau]$ has no zero divisors. But the right hand side is in $\Fm[\tau]$, because $\mu^{-1}\in\Fm$ and $\lambda^{-1}\xi_U$, $\phi_a \in R[\tau]$. This contradiction proves the claim.

Now expand $\phi'_a = \sum_i c_i\tau^i$ with $c_i\in R$. Then $\lambda^{-1}\phi'_a\lambda  = \sum_i c_i\lambda^{p^i-1}\tau^i$, and the claim asserts that  $c_i\lambda^{p^i-1}\in R$. Since $\phi'\mod\Fm$ is a Drinfeld module of rank $>0$ and $a\in A$ is non-constant, there exists an $i>0$ with $c_i \in R\setminus\Fm$. For this $i$ we then have $p^i-1>0$ and $\lambda^{p^i-1}\in R$, and hence $\lambda\in R$. But this contradicts the original choice of~$\lambda$, and the proposition is proved.
\end{Proof}

\begin{Prop}
\label{GDM:IsoExt}
In the situation of Proposition \ref{GDM:HomExt} we have:
\begin{itemize}
\item[(a)] If $\xi_U$ is an isogeny, then so is~$\xi$.
\item[(b)] If $\xi_U$ is an isomorphism, then so is~$\xi$.
\end{itemize}
\end{Prop}

\begin{Proof}
Assume that $\xi_U$ is an isogeny, and let $\eta$ denote the generic point of~$U$. Then there exists an isogeny in the reverse direction $\xi'_\eta\!: (E',\phi')|_\eta \to (E,\phi)|_\eta$ such that $\xi'_\eta\circ\xi|_\eta = \phi_a|_\eta$ for some non-zero element $a\in A$. This isogeny extends to a homomorphism over an open neighborhood of $\eta$ and then, by Proposition \ref{GDM:HomExt}, to a homomorphism $\xi'\!: (E',\phi') \to (E,\phi)$. By uniqueness, this extension still satisfies $\xi'\circ\xi = \phi_a$, and so the same relation holds in every fiber. It follows that $\xi$ is an isogeny in every fiber, proving~(a). 

If $\xi_U$ is an isomorphism, its inverse extends to a morphism $\xi^{-1}\!: (E',\phi') \to (E,\phi)$ by Proposition \ref{GDM:HomExt}, and by uniqueness both $\xi\circ\xi^{-1}$ and $\xi^{-1}\circ\xi$ are the identity. Thus $\xi$ and $\xi^{-1}$  are mutually inverse isomorphisms, proving~(b).
\end{Proof}

%As a direct consequence of Proposition \ref{GDM:IsoExt} (b) we deduce:
%\begin{Cor}
%\label{GDM:ObjExtUni}
%Assume that $S$ is normal integral with generic point~$\eta$. Then any generalized Drinfeld $A$-module over $\eta$ possesses at most one extension to a generalized Drinfeld $A$-module over $S$, up to a unique isomorphism extending the identity over~$\eta$.
%\end{Cor}

\medskip
For use in the next section we include the following notion:

\begin{Def}
\label{GDM:SepDef}
A generalized Drinfeld $A$-module $(E,\phi)$ over $S$ is \emph{weakly separating} if, for any Drinfeld $A$-module $(E',\phi')$ over any field $L$ containing $F$, at most finitely many fibers of $(E,\phi)$ over $L$-valued points of $S$ are isomorphic to $(E',\phi')$.
\end{Def}

\begin{Prop}
\label{GDM:Dim}
Let $(E,\phi)$ be a weakly separating generalized Drinfeld $A$-module over a scheme $S$ of finite type over~$F$. Then for any integer $r$ the set $S_r$ of points $s\in S$ where the fiber has rank $r_s\le r$ is Zariski closed of dimension $\le r-1$.
\end{Prop}

\begin{Proof}
By semicontinuity $S_r$ is Zariski closed, so it possesses a unique structure of reduced closed subscheme. Also by semicontinuity every irreducible component of $S_r$ contains an open dense subset $U$ over which the rank $r_s$ is constant, say equal to $r'\le r$. Then $(E,\phi)|U$ is a Drinfeld $A$-module of rank~$r'$. Thus for any non-zero proper ideal $N\subsetneqq A$, there exist a finite Galois covering $U'\onto U$ and a level $N$ structure on the  pullback of $(E,\phi)$ to~$U'$. By the modular interpretation this data corresponds to a morphism $f\!: U' \to M^{r'}_{A,K(N)}$. Since $(E,\phi)$ is weakly separating, so is its pullback to~$U'$. But by construction this pullback is also isomorphic to the pullback of the universal family under~$f$. Thus Definition \ref{GDM:SepDef} implies that the fibers of $f$ are finite. It follows that $\dim U \le \dim M^{r'}_{A,K(N)} = r'-1 \le r-1$. Therefore every irreducible component of $S_r$ has dimension $\le r-1$, as desired.
\end{Proof}

%%%%%%%%%%%%%%%%%%%%%%%%%%%%%%%%%%%%%%%%%%%%%%%%%%%%%%%%%%%%%%%%%%%%%%%%%%%%%%%%%%%%%%%%%%%%%

\section{Satake compactification}
\label{SC}

\begin{Def}
\label{SC:SatDef}
For any fine open compact subgroup $K \subset \GL_r(\hat A)$, an open embedding $M^r_{A,K} \into \OM^r_{A,K}$ with the properties
\begin{itemize}
\item[(a)] $\OM^r_{A,K}$ is a normal integral proper algebraic variety over~$F$, and
\item[(b)] the universal family on $M^r_{A,K}$ extends to a weakly separating generalized Drinfeld $A$-module $(\bar{E},\bar{\phi})$ over $\OM^r_{A,K}$,
\end{itemize}
is called a \emph{Satake compactification of $M^r_{A,K}$}. By abuse of terminology we call $(\bar{E},\bar{\phi})$ the \emph{universal family on $\OM^r_{A,K}$}.
\end{Def}

\begin{Thm}
\label{SC:SatExUni}
For every fine $K$ the variety $M^r_{A,K}$ possesses a projective Satake compactification. The Satake compactification and its universal family are unique up to unique isomorphism.
% . $\OM^r_{A,K}$ and $(\bar{E},\bar{\phi})$ are unique up to unique isomorphism.
\end{Thm}

Most of the proof resides in the following four lemmas:

%%%%%%%%%%%%%%%%%%%%%%%

\begin{Lem}
\label{SC:SatExUniLemUniq}
If $\OM^r_{A,K}$ and $(\bar{E},\bar{\phi})$ exist, they are unique up to unique isomorphism.
% The Satake compactification is unique up to unique isomorphism, if it exists.
\end{Lem}

\begin{Proof}
Abbreviate $\OM := \OM^r_{A,K}$ and let $\OM'$ be another Satake compactification of $M := M^r_{A,K}$ with universal family $(\bar{E}',\bar{\phi}')$. Let $\tilde M$ be the normalization of the Zariski closure of $M$ in $\OM\times_F\OM'$. Then the projections $\OM \stackrel{\pi}{\leftarrow} \tilde M \stackrel{\pi'}{\to} \OM'$ are proper and restrict to the identity on~$M$. By Propositions \ref{GDM:HomExt} and \ref{GDM:IsoExt}~(b) the identity on the universal family on $M$ extends to an isomorphism $\pi^*(\bar{E},\bar{\phi}) \cong \pi^{\prime*}(\bar{E}',\bar{\phi}')$. Thus for any geometric point $x\in\OM(L)$ over an algebraically closed field~$L$, the isomorphism class of $(\bar{E}',\bar{\phi}')$ is constant over the set of points $\pi'(\pi^{-1}(x)) \allowbreak \subset\nobreak \OM'$. Since $(\bar{E}',\bar{\phi}')$ is weakly separating, it follows that $\pi'(\pi^{-1}(x))$ is finite. By the construction of $\tilde M$ this implies that the morphism $\pi$ is quasi-finite and hence finite.
% Suppose that $\pi$ is not finite. Then there is a closed irreducible subvariety $X\subset \tilde M$ which is mapped to a point under~$\pi$. The morphism $\pi'$ is then non-constant on $X$ by the construction of~$\tilde M$. Thus if $\pi$ had a fiber of dimension $>0$, ...
As $\OM$ is already normal and $\pi$ is birational, we deduce that $\pi$ is an isomorphism. In the same way one proves that $\pi'$ is an isomorphism. Thus $(\OM,\bar{E},\bar{\phi})$ is isomorphic to $(\OM',\bar{E}',\bar{\phi}')$, and clearly the isomorphism extending the identity is unique.
\end{Proof}

%%%%%%%%%%%%%%%%%%%%%%%

\begin{Lem}
\label{SC:SatExUniLemPush}
For any two fine open compact subgroups of the form $K(N)\subset K$, if $M^r_{A,K(N)}$ possesses a projective Satake compactification $\OM^r_{A,K(N)}$, then $M^r_{A,K}$ possesses the projective Satake compactification $\OM^r_{A,K} := \OM^r_{A,K(N)}/K$.
\end{Lem}

\begin{Proof}
Recall from Proposition \ref{DMV:FineProp} that $M^r_{A,K}$ is the quotient of $M^r_{A,K(N)}$ under a free action of the finite group $K/K(N)$. This group also acts on the universal family by the isomorphisms (\ref{MI:R3}); hence by Lemma \ref{SC:SatExUniLemUniq} the action extends to $\OM^r_{A,K(N)}$ and its universal family $(\bar E',\bar\phi')$. As $\OM^r_{A,K(N)}$ is projective, the desired quotient $\OM^r_{A,K}$ exists and is a normal integral projective algebraic variety over~$F$ containing $M^r_{A,K}$ as an open subvariety. 

Since $\bar E'$ is quasi-projective, we can also form the quotient $\bar E := \bar E'/K$ as an algebraic variety. We claim that $\bar E$ is a line bundle on $\OM^r_{A,K}$ whose pullback to $\OM^r_{A,K(N)}$ is~$\bar E'$. Granting this for the moment, the equi\-va\-ri\-ance implies that the morphisms $\bar\phi'_a\!: \bar E'\to\bar E'$ are the pullbacks of morphisms $\bar\phi_a\!: \bar E\to\bar E$ for all $a\in A$. Thus $(\bar E',\bar\phi')$ is the pullback of $(\bar E,\bar\phi)$, and so the latter is a generalized Drinfeld $A$-module over $\OM^r_{A,K}$. Since $(\bar E',\bar\phi')$ is already weakly separating, the same follows also for $(\bar{E},\bar{\phi})$. Moreover $(\bar E,\bar\phi)$ restricts to the universal family on $M^r_{A,K}$ by Proposition~\ref{DMV:FineProp}. Thus $\OM^r_{A,K}$ satisfies all the conditions of \ref{SC:SatDef}, and the proposition follows.

To prove the claim, standard descent theory asserts that it suffices to work \'etale locally on $\OM^r_{A,K}$. 
% By standard descent theory it suffices to prove the claim \'etale locally on $\OM^r_{A,K}$. 
Thus we may replace $K/K(N)$ by the stabilizer $\Delta_x \subset K/K(N)$ of a geometric point $x$ of $\OM^r_{A,K(N)}$, and $\OM^r_{A,K(N)}$ by a $\Delta_x$-invariant \'etale neighborhood $U_x$ of~$x$. By assumption the fiber $(\bar E'_x,\bar\phi'_x)$ over $x$ is a Drinfeld module of some rank $r_x>0$. Thus its $N$-torsion points form a free $A/N$-module of rank~$r_x$. Moreover $\bar\phi'[N]$ is an \'etale (though not necessarily finite) group scheme over $\OM^r_{A,K(N)}$, because for any non-zero $a\in A$, the coefficient of $\tau^0=1$ of $\bar\phi'_a$ is the image of $a\in F^\times$ by \ref{GDM:GenDrinDef} (a) and hence non-zero everywhere. 
% (But $\bar\phi'[N]$ is not finite where the rank decreases.)
Thus every $N$-torsion point of $\bar\phi'_x$ extends to a section of $\bar\phi'[N]$ over an \'etale neighborhood of~$x$. Let $\lambda\!:$ $\underline{(N^{-1}/A)^r} \stackrel\sim\longto \bar\phi'[N]|_{M^r_{A,K(N)}}$ be the given level $N$ structure over $M^r_{A,K(N)}$, and let $W$ denote the set of $w\in (N^{-1}/A)^r$ for which $\lambda(w)$ extends to a section of $\bar\phi'[N]$ over a neighborhood of~$x$. Then the preceding remarks show that $W$ is a free $A/N$-submodule of rank~$r_x$.

Since $K$ is fine, by Definition \ref{DMV:FineDef} the image of $K$ in $\GL_r(A/\Fp)$ is unipotent for some maximal ideal $\Fp\subset A$. As $K(N)\subset K$, we must have $N\subset\Fp$. Thus $K$ and hence $\Delta_x$ act unipotently on $(\Fp^{-1}/A)^r$. On the other hand $W\cap(\Fp^{-1}/A)^r$ is a free $A/\Fp$-module of rank $r_x>0$ and by construction invariant under~$\Delta_x$. We can therefore find a non-zero element $w\in W\cap(\Fp^{-1}/A)^r$ which is fixed by~$\Delta_x$. 

By construction $\lambda(w)$ extends to a section $\bar\lambda(w)$ of $\bar\phi'[N]$ over a neighborhood $U_x$ of~$x$, which we can also take $\Delta_x$-invariant. Since $\lambda(w)$ is non-zero and $\bar\phi'[N]$ is \'etale, the extension $\bar\lambda(w)$ is non-zero everywhere. Thus it defines a $\Delta_x$-invariant trivialization $\BG_{a,U_x} \stackrel{\sim}{\to} \bar E'|_{U_x}$. It follows that this trivialized line bundle is the  pullback of a trivialized line bundle on $U_x/\Delta_x$. The latter can of course be constructed as $\BG_{a,U_x/\Delta_x} = \BG_{a,U_x}/\Delta_x \cong (\bar E'|_{U_x})/\Delta_x$ and therefore has the desired properties.
\end{Proof}

%%%%%%%%%%%%%%%%%%%%%%%

\begin{Lem}
% [Pullback]
\label{SC:SatExUniLemPull}
For any morphism $I_b\!: M^{r'}_{A',K(N')} \to M^r_{A,K(N)}$ as in (\ref{MI:Theta1}), if a projective Satake compactification exists for $M^r_{A,K(N)}$, then one exists for $M^{r'}_{A',K(N')}$.
\end{Lem}

\begin{Proof}
Let $\OM^r_{A,K(N)}$ be a projective Satake compactification of $M^r_{A,K(N)}$, and define $\OM^{r'}_{A',K(N')}$ as the normalization of $\OM^r_{A,K(N)}$ in the function field of $M^{r'}_{A',K(N')}$. Since the morphism $I_b$ is finite and the scheme $M^{r'}_{A',K(N')}$ is normal, we obtain a commutative diagram
$$\xymatrix{
M^{r'}_{A',K(N')} \ar@{^{ (}->}[r] \ar[d]^{I_b} & \OM^{r'}_{A',K(N')} \ar[d]^{\bar I_b} \\
M^{r }_{A ,K(N )} \ar@{^{ (}->}[r]              & \OM^{r }_{A ,K(N )}\\}$$
where the horizontal arrows are open embeddings and the morphism $\bar I_b$ is finite. Then $\OM^{r'}_{A',K(N')}$ is projective; we will show that it is a Satake compactification of $M^{r'}_{A',K(N')}$.

Let $(\bar{E},\bar{\phi})$ denote the universal family on $\OM^r_{A,K(N)}$ and $(E',\phi')$ the universal family on $M^{r'}_{A',K(N')}$. Then $(\bar E',\tilde\phi) := \bar I_b^*(\bar{E},\bar{\phi})$ is a generalized Drinfeld $A$-module over $\OM^{r'}_{A',K(N')}$, whose restriction to $M^{r'}_{A',K(N')}$ is isomorphic to $(E',\phi'|A)$ by (\ref{MI:Theta2}). Since $A'$ is commutative, we can view the endomorphism $\phi'_a\!:E'\to E'$ associated to any $a\in A$ as an endomorphism of the Drinfeld $A$-module $(E',\phi'|A)$. By Proposition \ref{GDM:HomExt} it therefore extends to a unique endomorphism of $(\bar E',\tilde\phi)$. Again by uniqueness, this defines an algebra homomorphism $A\to \End(\bar E')$ and thus a generalized Drinfeld $A$-module $(\bar{E}',\bar{\phi}')$ extending $(E',\phi')$ such that $\bar{\phi}'|A = \tilde\phi$. Since $(\bar{E},\bar{\phi})$ is weakly separating and $\bar I_b$ is finite, this implies that $(\bar{E}',\bar{\phi}'|A)$ is weakly separating. Therefore $(\bar{E}',\bar{\phi}')$ is weakly separating, and so $\OM^{r'}_{A',K'}$ and $(\bar{E}',\bar{\phi}')$ satisfy all the desired conditions.
\end{Proof}

%%%%%%%%%%%%%%%%%%%%%%%

\medskip
The following special case will be proved in Section \ref{Fqtt}: see Theorem \ref{Fqtt:Sat}.

\begin{Lem}
\label{SC:SatExUniLemBase}
A projective Satake compactification exists for $M^r_{\BF_q[t],K(t)}$ for any $r\ge1$ and any finite extension $\BF_q$ of~$\BF_p$.
\end{Lem}

%%%%%%%%%%%%%%%%%%%%%%%

\begin{Proofof}{Theorem \ref{SC:SatExUni}}
The uniqueness part is contained in Lemma \ref{SC:SatExUniLemUniq}. To construct a projective Satake compactification of $M^r_{A,K}$ take a principal congruence subgroup $K(N) \subset\nobreak K$. Choose a non-zero element $t\in N$; then after shrinking $N$ we may assume that $N=(t)$. By Lemma \ref{SC:SatExUniLemPush} it suffices to show that a projective Satake compactification exists for $M^r_{A,K(N)}$. For this set $A' := \BF_p[t] \subset A$ and $N' := (t) \subsetneqq A'$, and choose an $\hat A'$-linear isomorphism $b\!: \hat A^{\prime\kern1pt r'} \stackrel{\sim}{\to} \hat A^r$. Then by Lemma \ref{SC:SatExUniLemPull} (with primed and unprimed letters interchanged) it suffices to show that a projective Satake compactification exists for $M^{r'}_{\BF_p[t],K(t)}$. But this is guaranteed by Lemma \ref{SC:SatExUniLemBase}.
\end{Proofof}

%%%%%%%%%%%%%%%%%%%%%%%

\medskip
Now consider an open compact subgroup $K\subset \GL_r(\hat A)$ which is not fine. Then we cannot characterize a compactification of $M^r_{A,K}$ in terms of a universal family as in Definition~\ref{SC:SatDef}. Instead we choose any $K(N)\subset K$. The uniqueness in Theorem~\ref{SC:SatExUni} implies that the action of $K/K(N)$ on $M^r_{A,K(N)}$ and its universal family extends to an action on $\OM^r_{A,K(N)}$. Since the latter is projective, we can form the quotient
\addtocounter{Thm}{1}
\begin{equation}
\label{SC:SatDefNotFine}
\OM^r_{A,K} := \OM^r_{A,K(N)}/K.
\end{equation}
This is a normal integral projective algebraic variety over~$F$, and by (\ref{DMV:DMVDef}) it contains the open subvariety $M^r_{A,K}$. We call it the \emph{Satake compactification of $M^r_{A,K}$} in this case. It is independent of the choice of~$N$, because 
% the isomorphism (\ref{DMV:Quot}) extends to the Satake compactifications
$$\OM^r_{A,K(N)} \cong \OM^r_{A,K(N')}/K(N)$$
for any $N'\subset N$, which is a special case of the following fact:

\begin{Prop}
\label{SC:SatQuot}
For any open compact subgroup $K\subset\GL_r(\hat A)$ and any open normal subgroup $K'\triangleleft K$ we have 
$$\OM^r_{A,K} \cong \OM^r_{A,K'}/K.$$
\end{Prop}

\begin{Proof}
If $K$ and hence $K'$ are fine, take any $K(N) \subset K'$. Then Lemma \ref{SC:SatExUniLemPush} implies that $\OM^r_{A,K} \cong \OM^r_{A,K(N)}/K = \bigl(\OM^r_{A,K(N)}/K'\bigr)/K \cong \OM^r_{A,K'}/K$, as desired. The general case follows in the same way using the definition (\ref{SC:SatDefNotFine}).
\end{Proof}

%%%%%%%%%%%%%%%%%%%%%%%

\begin{Rem}
\label{SC:Kap}
\rm Kapranov \cite[Thm.$\,$1.1]{Kapranov} already constructed a projective compactification of $M^r_{A,K(N)}$ in the case $A=\BF_q[t]$ and proved that its boundary is a union of finitely many moduli spaces of the form $M^{r'}_{A,K(N)}$ for $1\leq r'<r$. 
% and suitable $K'\subset\GL_{r'}(\Ahat)$.
We are convinced that this is a Satake compactification in our sense, but the proof is necessarily based on Kapranov's analytic construction and thus outside the scope of the present article.
\end{Rem}

In the general case, too, we expect $\OM^r_{A,K}$ to be stratified by finitely many moduli spaces of the form $M^{r'}_{A,K'}$ for different $r'$ and~$K'$. This would directly imply the following result, which we prove independently:
% because we need it for the application to modular forms.

\begin{Prop}
\label{SC:Div}
If $K$ is fine, the fiber of the universal family over the generic point of any irreducible component of $\OM^r_{A,K} \setminus M^r_{A,K}$ is a Drinfeld $A$-module of rank ${r-1}$.
\end{Prop}

\begin{Proof}
By (\ref{SC:SatDefNotFine}) it suffices to prove this when $K=K(N)$. Let $X$ be the irreducible component in question and $x$ its generic point. Let $r_x$ be the rank of the universal family above~$x$. Then by semicontinuity the rank is $\le r_x$ over all of~$X$, and Proposition \ref{GDM:Dim} implies that $\dim X \le r_x-1$. On the other hand, the fact that $M^r_{A,K(N)}$ is affine of dimension $r-1$ (see \cite[\S5]{Drinfeld1}, \cite[Thm.$\,$1.8]{GekelerOhio}) implies that $\dim X = r-2$. Together we find that $r-1\le r_x$. 

Suppose that $r_x>r-1$. Then by semicontinuity we have $r_x=r$, and so the universal family $(\bar E,\bar\phi)$ is a Drinfeld $A$-module of rank $r$ over a neighborhood $U$ of~$x$. Thus the scheme of $N$-torsion points $\bar\phi[N]$ is finite \'etale over~$U$. By the valuative criterion for properness, the level $N$ structure over $M^r_{A,K(N)}$ therefore extends to a level $N$ structure over the local ring at~$x$. By the modular interpretation of $M^r_{A,K(N)}$ this means that $x$ is really a point of $M^r_{A,K(N)}$, and not its complement, contrary to the assumption. This proves that $r_x=r-1$, as desired.
\end{Proof}

%%%%%%%%%%%%%%%%%%%%%%%

\medskip
The morphisms and isogenies from Section \ref{MI} extend to the Satake compactifications in the following way.

\begin{Prop}
\label{SC:RgSat}
Consider elements $g$, $g'\in\GL_r(\BA_F^f)$ and open compact subgroups $K$, $K'$, $K''\subset\GL_r(\hat A)$ such that $gK'g^{-1}\subset K$ and $g'K''g^{\prime-1}\subset K'$. 
\begin{itemize}
\item[(a)] The morphism $J_g$ from Proposition \ref{MI:FineProp1} extends to a unique finite morphism
$$\bar J_g\!:\ \OM^r_{A,K'} \longto \OM^r_{A,K}.$$
\item[(b)] The extensions in (a) satisfy\ \ $\bar J_g\circ\bar J_{g'} =\bar J_{gg'}$.
\end{itemize}
Furthermore assume that $g$ and $g'$ have coefficients in~$\hat A$ and that $K$, $K'$, $K''$ are fine. Let $(\bar E,\bar\phi)$ and $(\bar E',\bar\phi')$ denote the respective universal families on $\OM^r_{A,K}$ and $\OM^r_{A,K'}$. 
\begin{itemize}
\item[(c)] The isogeny $\xi_g$ from Proposition \ref{MI:FineProp2} extends to a unique isogeny 
$$\bar\xi_g\!:\ (\bar E',\bar\phi') \longto \bar J_g^*(\bar E,\bar\phi).$$
\item[(d)] The extensions in (c) satisfy\ \ $(\bar J_{g'}^*\bar\xi_g)\circ \bar\xi_{g'} = \bar\xi_{gg'}$.
\end{itemize}
\end{Prop}

\begin{Proof}
By the construction of $J_g$ in the proof of Proposition \ref{MI:FineProp1}, it suffices to prove all this when $g$ and $g'$ have coefficients in~$\hat A$.

Assume first that $K$ and $K'$ are fine. Abbreviate $\OM := \OM^r_{A,K}$ and $\OM' := \OM^r_{A,K'}$, let $\tilde M'$ be the normalization of the closure of $\graph(J_g) \subset \OM'\times_F\OM$, and consider the proper projections $\smash{\OM' \stackrel{\pi'}{\leftarrow} \tilde M' \stackrel{\pi}{\to} \OM}$. Then $\pi^{\prime*}(\bar{E}',\bar{\phi}')$ and $\pi^*(\bar{E},\bar{\phi})$ are generalized Drinfeld $A$-modules over~$\tilde M'$, and by Propositions \ref{GDM:HomExt} and \ref{GDM:IsoExt}~(a) the isogeny $\xi_g$ extends to an isogeny $\tilde\xi_g\!:\ \pi^{\prime*}(\bar E',\bar\phi') \longto \pi^*(\bar E,\bar\phi)$. 
Since the degree of $\tilde\xi_g$ in a fiber defines a constructible function on~$\tilde M'$, it is bounded. Let $c$ be an upper bound for it.
% Let $c$ be an upper bound for the degree of $\tilde\xi_g$ in all fibers.

For any geometric point $x\in\OM(L)$ over an algebraically closed field~$L$, the isomorphism class of $\pi^*(\bar E,\bar\phi)$ is constant over $\pi^{-1}(x)$. Since any Drinfeld module over a field possesses only finitely many isogenies up to isomorphism into it, or out of it, of degree $\le c$, the fibers of $\pi^{\prime*}(\bar E',\bar\phi')$ over $\pi^{-1}(x)$ form only finitely many isomorphism classes. Thus the fibers of $(\bar{E}',\bar{\phi}')$ over $\pi'(\pi^{-1}(x)) \allowbreak \subset\nobreak \OM'$ fall into only finitely many isomorphism classes. Since $(\bar{E}',\bar{\phi}')$ is weakly separating, it follows that $\pi'(\pi^{-1}(x))$ is finite. By the construction of $\tilde M$ this implies that the morphism $\pi$ is quasi-finite and hence finite. In the same way one proves that $\pi'$ is finite. 

But as $\OM'$ is already normal and $\pi'$ is birational, this implies that $\pi'$ is an isomorphism. Thus $\pi\circ\pi^{\prime-1}$ is the desired extension in (a) and $\pi'_*\tilde\xi_g$ is the desired isogeny in~(c).
In particular this proves~(c), and it proves (a) when the subgroups are sufficiently small. The general case of (a) now results from Proposition~\ref{SC:SatQuot}. Part (b) follows from (\ref{MI:R2xa}) and the fact that $\OM^r_{A,K''}$ is normal integral. Part~(d) follows from the last statement in Proposition \ref{MI:FineProp2} and the uniqueness part of Proposition~\ref{GDM:HomExt}. 
\end{Proof}

%%%%%%%%%%%%%%%%%%%%%%%

\begin{Prop}
\label{SC:TeSat}
\begin{itemize}
\item[(a)] The morphism $I_b$ of (\ref{MI:Theta3}) extends to a unique finite morphism 
$$\bar I_b\!:\ \OM^{r'}_{A',K'} \longto \OM^r_{A,K}.$$
\item[(b)] If $K$ and $K'$ are fine and $(\bar E,\bar\phi)$ and $(\bar E',\bar\phi')$ denote the respective universal families on $\OM^r_{A,K}$ and $\OM^{r'}_{A',K'}$, there is a natural isomorphism
$$(\bar E',\bar\phi'|A)\ \cong\ \bar I_b^*(\bar E,\bar\phi).$$
\end{itemize}
\end{Prop}

\begin{Proof}
We may assume that $K$ and $K'$ are fine; the general case of (a) then follows using Proposition~\ref{SC:SatQuot}. Abbreviate $\OM := \OM^r_{A,K}$ and $\OM' := \OM^{r'}_{A',K'}$, let $\tilde M'$ be the normalization of the closure of $\graph(I_b) \subset \OM'\times_F\OM$, and consider the proper projections $\smash{\OM' \stackrel{\pi'}{\leftarrow} \tilde M' \stackrel{\pi}{\to} \OM}$. By Propositions \ref{GDM:HomExt} and \ref{GDM:IsoExt}~(b) the isomorphism (\ref{MI:Theta2}) yields an isomorphism of generalized Drinfeld $A$-modules $\pi^{\prime*}(\bar{E}',\bar{\phi}') \cong \pi^*(\bar{E},\bar{\phi}|A') $ over~$\tilde M'$.

For any geometric point $x\in\OM(L)$ over an algebraically closed field~$L$, it follows that the isomorphism class of $(\bar{E}',\bar{\phi}'|A)$ is constant over $\pi'(\pi^{-1}(x)) \subset \OM'$. Since a Drinfeld $A$-module over a field possesses at most finitely many extensions to a Drinfeld $A'$-module, the fibers of $(\bar{E}',\bar{\phi}')$ over $\pi'(\pi^{-1}(x))$ fall into only finitely many isomorphism classes. As $(\bar{E}',\bar{\phi}')$ is weakly separating, we deduce that $\pi'(\pi^{-1}(x))$ is finite. By the construction of $\tilde M'$ this implies that $\pi$ is finite. In the same way one proves that $\pi'$ is finite.

But $\OM'$ is already normal and $\pi'$ is birational; hence $\pi'$ is an isomorphism, and the proposition follows.
\end{Proof}

%%%%%%%%%%%%%%%%%%%%%%%%%%%%%%%%%%%%%%%%%%%%%%%%%%%%%%%%%%%%%%%%%%%%%%%%%%%%%%%%%%%%%%%%%%%%%

\section{Modular forms}
\label{MF}

First we assume that $K$ is fine. Let $(\bar E,\bar\phi)$ be the universal family over the Satake compactification $\OM^r_{A,K}$. Let $\smash{\CL_{\OM^r_{A,K}}}$ denote the dual of the 
relative Lie algebra of $\bar E \onto \OM^r_{A,K}$, 
% sheaf of sections of~$\bar E$,
which is an invertible sheaf on $\OM^r_{A,K}$.
% ; we often abbreviate it by~$\CL$.

\begin{Lem}
\label{MF:JgIe}
In the situation of Proposition \ref{SC:RgSat} (c--d) and \ref{SC:TeSat} we have:
\begin{itemize}
\item[(a)]
The derivative of $\bar\xi_g$ induces a natural isomorphism 
$$\rho_g\!:\ \ \bar J_g^*\CL_{\OM^r_{A,K}} \stackrel{\sim}{\longto} \CL_{\OM^r_{A,K'}}.$$
\item[(b)]
The isomorphisms in (a) satisfy\ \ $\rho_{g'} \circ (\bar J_{g'}^*\rho_g)  = \rho_{gg'}$.
\item[(c)]
There is a natural isomorphism\ \ $\bar I_b^*\CL_{\OM^r_{A,K}} \cong \CL_{\OM^{r'}_{A',K'}}$.
\end{itemize}
\end{Lem}

\begin{Proof}
Every isogeny of Drinfeld modules of rank $>0$ over an extension of $F$ is separable, i.e., its derivative is non-zero (see \cite[Prop.$\,$4.7.10]{GossBS}). Thus the derivative of $\bar\xi_g$ is an isomorphism of Lie algebras in all fibers, and its dual is the isomorphism~$\rho_g$. This proves (a), and the cocycle relation (b) follows at once from \ref{SC:RgSat}~(d). Part (c) is proved in the same way as~(a).
\end{Proof}

\begin{Rem}
\label{MF:WhyLie}
\rm Part (a) of the above lemma explains why the natural definition of modular forms requires the relative Lie algebra of the line bundle $\bar E \onto \OM^r_{A,K}$ and not the line bundle itself:
% (although they are isomorphic): 
namely, because $d\bar\xi_g$ is an isomorphism of invertible sheaves whereas $\bar\xi_g$ is in general not an isomorphism of line bundles.
\end{Rem}

\begin{Thm}
\label{MF:Ample}
$\CL_{\OM^r_{A,K}}$ is ample.
\end{Thm}

\begin{Proof}
Take any principal congruence subgroup $K(N) \subset K$. Then Proposition \ref{SC:SatQuot} and Lemma \ref{MF:JgIe}~(a) imply that the assertion holds for $\OM^r_{A,K}$ if it holds for $\OM^r_{A,K(N)}$. Next consider the morphism $\bar I_b\!: \OM^r_{A,K(N)} \to \OM^{r'}_{\BF_p[t],K(t)}$ as in the proof of Theorem~\ref{SC:SatExUni}. Since that morphism is finite by Proposition~\ref{SC:TeSat}, Lemma \ref{MF:JgIe}~(c) reduces the assertion to the case $\OM^{r'}_{\BF_p[t],K(t)}$. In that case the assertion is part of Theorem \ref{Fqtt:Sat} below.
\end{Proof}

%%%%%%%%%%%%%%%%%%%%%%%

\medskip
In the special case $K=K'=K(N)$ and $g\in\GL_r(\hat A)$ Lemma \ref{MF:JgIe} (a--b) shows that the action of $\GL_r(\hat A)$ on $\OM^r_{A,K(N)}$ lifts to a covering action on $\CL_{\OM^r_{A,K(N)}}$. Thus it induces an action on global sections, so that the following definition makes sense:

\begin{Def}
\label{MF:MFDef}
For any integer $k$ we call
$$\qquad\ \CM_k(M^r_{A,K}) \ :=\ 
\begin{cases}
H^0\bigl(\OM^r_{A,K},\CL^k_{\OM^r_{A,K}}\bigr) & \hbox{if $K$ is fine,}\\
H^0\bigl(\OM^r_{A,K(N)},\CL^k_{\OM^r_{A,K(N)}}\bigr)^K & 
\hbox{if $K$ is not fine, but $K(N)\subset K$ is.}
\end{cases}$$
the \emph{space of (algebraic) modular forms of weight $k$ on $M^r_{A,K}$}, and
$$\CR(M^r_{A,K}) \ :=\ \bigoplus_{k\ge0}\; \CM_k(M^r_{A,K})$$
the \emph{ring of (algebraic) modular forms on $M^r_{A,K}$}.
\end{Def}

When $K$ is not fine, this definition is independent of~$N$, because $\CM_k(M^r_{A,K(N)}) = \CM_k(M^r_{A,K(N')})^{K(N)}$ for any $N'\subset N$, which is a special case of the following fact:

\begin{Prop}
\label{MF:Quot}
For any open compact subgroup $K\subset\GL_r(\hat A)$ and any open normal subgroup $K'\triangleleft K$ we have a natural isomorphism
$$\CM_k(M^r_{A,K})\ \cong\ \CM_k(M^r_{A,K'})^K.$$
\end{Prop}

\begin{Proof}
Proposition \ref{SC:SatQuot} and Lemma~\ref{MF:JgIe} imply this when $K$ and $K'$ are fine. The general case follows from this and Definition~\ref{MF:MFDef}.
\end{Proof}

%%%%%%%%%%%%%%%%%%%%%%%

\begin{Thm}
\label{MF:Ring}
The ring $\CR(M^r_{A,K})$ is a normal integral domain that is finitely generated as an $F$-algebra, and
$$\OM^r_{A,K} \cong \Proj \CR(M^r_{A,K}).$$
\end{Thm}

\begin{Proof}
When $K$ is fine, this follows from Theorem \ref{MF:FinGen} below, which for lack of a suitable reference we include with a proof. The general case follows from this and Noether's theorem \cite[Thm.$\,$2.3.1]{LarrySmith} that for any finite group acting on a finitely generated algebra, the subring of invariants is again finitely generated. 
\end{Proof}

\begin{Thm}
\label{MF:FinGen}
For any normal integral projective algebraic variety $X$ over a field~$F$, the ring $R$ of sections in all powers of an ample invertible sheaf $\CL$ is a finitely generated normal integral domain, and $X=\Proj R$. 
\end{Thm}

\begin{Proof}
For each $i\ge0$, the space $R_i := H^0(X,\CL^i)$ is finite dimensional over $F$ by the coherence theorem \cite[Thm.$\,$5.19]{Hartshorne}. Fix an integer $n\ge0$ such that $\CL^n$ is very ample. Then we have a natural short exact sequence $0\to \CF\to R_n\otimes_F\CO_X\to\CL^n\to0$ with a coherent sheaf $\CF$ on~$X$. Since $\CL$ is ample, there exists an integer $i_0$ such that $H^1(X,\CF\otimes\CL^i)=0$ for all $i\ge i_0$. For these $i$ the long exact cohomology sequence associated to the short exact sequence $0\to \CF\otimes\CL^i\to R_n\otimes_F\CL^i\to\CL^{n+i}\to0$ implies that the multiplication map $R_n\otimes R_i \cong H^0(X,R_n\otimes\CL^i) \to H^0(X,\CL^{n+i}) = R_{n+i}$ is surjective. It follows that $R = \bigoplus_{i\ge0}R_i$ is generated as an $F$-algebra by $R_n$ together with $R_i$ for all $i\le i_0$; hence it is finitely generated.

By the argument in \cite[Ex.$\,$5.14$\,$(a)]{Hartshorne}, which requires only an ample invertible sheaf, not necessarily a very ample one, the ring $R$ is a normal integral domain.
%Next, since $X$ is integral, any non-zero element of $R_i$ is invertible at the generic point of~$X$. Thus the product of any two non-zero homogeneous elements of $R$ is non-zero. Expanding any two non-zero elements $x$, $y\in R$ into their homogeneous components, we find that the product of their respective highest non-zero components is non-zero; hence $xy$ is non-zero, and so $R$ is an integral domain.
%
%Next, the normalization of any graded integral domain is itself graded by \cite[Thm.$\,$2.3.2]{HunekeSwanson}. To prove that $R$ is normal it thus suffices to show that any quotient $x/y$ of non-zero homogeneous elements of~$R$, which is integral over $R$, lies in~$R$. For this suppose that $x\in R_i$ and $y\in R_j$. . . . 
\end{Proof}

%%%%%%%%%%%%%%%%%%%%%%%%%%%%%%%%%%%%%%%%%%%%%%%%%%%%%%%%%%%%%%%%%%%%%%%%%%%%%%%%%%%%%%%%%%%%%

\section{Hecke operators}
\label{HO}

Consider an element $g\in\GL_r(\BA_F^f)$ with coefficients in~$\hat A$, and open compact subgroups $K$, $K'\subset\GL_r(\hat A)$ such that $gK'g^{-1}\subset K$. For any such data we want to construct a \emph{pullback map}
\addtocounter{Thm}{1}
\begin{equation}
\label{MF:Pull}
J_g^*\!:\ \ \CM_k(M^r_{A,K})\ \longto\ \CM_k(M^r_{A,K'}).
\end{equation}
If $K$ and $K'$ are fine, this map is defined by composing the pullback of sections with the isomorphism $\rho_g$ from Lemma \ref{MF:JgIe}~(a). In the special case $K=K'=K(N)$ and $g\in\GL_r(\hat A)$ this map yields the group action used in Section~\ref{MF}. In the case $g=1$ it is the inclusion $\CM_k(M^r_{A,K}) \into \CM_k(M^r_{A,K'})$ from Proposition~\ref{MF:Quot}.

\begin{Prop}
\label{MF:PullGen}
There is a unique way of defining $J_g^*$ in the remaining cases such that $J_{g'}^* \circ J_g^*  = J_{gg'}^*$ whenever the formula makes sense.
%For any further element $g'\in\GL_r(\BA_F^f)$ with coefficients in $\hat A$ and any open compact subgroup $K''\subset\GL_r(\hat A)$ with $g'K''g^{\prime-1}\subset K'$, these pullback maps satisfy\ \ $J_{g'}^* \circ J_g^*  = J_{gg'}^*$.
\end{Prop}

\begin{Proof}
(Sketch) When all open compact subgroups involved are fine, the formula $J_{g'}^* \circ J_g^*  = J_{gg'}^*$ follows from Proposition \ref{SC:RgSat} (b) and Lemma \ref{MF:JgIe}~(b). When $K$ or $K'$ are not fine, choose fine open normal subgroups $L\triangleleft K$ and $L'\triangleleft K'$ such that $gL'g^{-1}\subset L$. Then the formula just proved together with some calculation implies that $J_g^*\!: \CM_k(M^r_{A,L}) \to \CM_k(M^r_{A,L'})$ sends $K$-invariants to $K'$-invariants. Using Proposition~\ref{MF:Quot} this yields the desired map (\ref{MF:Pull}). Further direct calculation then shows the formula $J_{g'}^* \circ J_g^*  = J_{gg'}^*$ in general.
\end{Proof}

\medskip
In the same way one can define a natural \emph{restriction map}
\addtocounter{Thm}{1}
\begin{equation}
\label{MF:Res}
\CM_k(M^r_{A,K})\ \longto\ \CM_k(M^{r'}_{A',K'})
\end{equation}
in the situation of (\ref{MI:Theta3}); namely by composing the pullback of sections with the isomorphism from  Lemma \ref{MF:JgIe}~(c) if $K$ and $K'$ are fine, and by taking invariants as in the proof of Proposition~\ref{MF:PullGen} in the general case. There is also a certain compatibility relation between the restriction maps (\ref{MF:Res}) and pullback maps (\ref{MF:Pull}) coming from elements of $\GL_{r'}(\BA_F^f)$, and a straightforward associativity relation for composites of restriction maps (\ref{MF:Res}), which the reader may write out and verify for him- or herself.

%%%%%%%%%%%%%%%%%%%%%%%

\medskip
Now we return to the situation of (\ref{MF:Pull}) and construct a map in the other direction. {}From (\ref{MI:R1}) we can see that $J_g\!: M^r_{A,K'} \to M^r_{A,gK'g^{-1}}$ is always an isomorphism. As the Satake compactification is normal integral, Proposition \ref{SC:RgSat} (a) implies that $\bar J_g\!: \OM^r_{A,K'} \to \OM^r_{A,gK'g^{-1}}$ is an isomorphism, too. It follows that $J_g^*\!: \CM_k(M^r_{A,gK'g^{-1}}) \to \CM_k(M^r_{A,K'})$ is an isomorphism. Choose an open normal subgroup $L\triangleleft K$ such that $L \subset gK'g^{-1}$. Then we have the commutative diagram
\vskip-10pt
$$\xymatrix@R-5pt@C+20pt{
\CM_k(M^r_{A,K})   \ar@{^{ (}->}[r]^-{J_1^*} \ar@{=}[d]^\wr_{\ref{MF:Quot}} \ar@/^5ex/[rr]^-{J_g^*} & 
\CM_k(M^r_{A,gK'g^{-1}})     \ar[r]^-{J_g^*} \ar@{=}[d]^\wr_{\ref{MF:Quot}} &
\CM_k(M^r_{A,K'}) \ar@{=}[d]^\wr_{\ref{MF:Quot}} \\
\CM_k(M^r_{A,L})^K \ar@{^{ (}->}[r]^-{J_1^*} &
\CM_k(M^r_{A,L})^{gK'g^{-1}} \ar[r]^-{J_g^*} \ar@{..>}@/^3ex/[l]^{\rm trace} &
\CM_k(M^r_{A,g^{-1}Lg})^{K'} \rlap{,}\\}$$
disregarding the dotted arrow that has not yet been defined. By the preceding remarks the two horizontal morphisms on the right hand side are isomorphisms. We can now define the dotted arrow as
\addtocounter{Thm}{1}
\begin{equation}
\label{MF:Trace}
f\longmapsto \mathop{\rm trace}(f) := \textstyle \sum_h J_h^*f,
\end{equation}
where $h$ runs through a set of representatives of the quotient $gK'g^{-1} \backslash K$.
%\begin{Def}
%\label{MF:Push}
%The \emph{pushforward map}
%$$J_{g\kern1pt*}\!:\ \ \CM_k(M^r_{A,K'})\ \longto\ \CM_k(M^r_{A,K})$$
%is the composite of the trace map with the isomorphisms \ref{MF:Quot} and $J_g^*$ in the above diagram.
%\end{Def}
The composite of this trace map with the isomorphisms \ref{MF:Quot} and $J_g^*$ in the above diagram then defines the \emph{pushforward map}
\addtocounter{Thm}{1}
\begin{equation}
\label{MF:Push}
J_{g\kern1pt*}\!:\ \ \CM_k(M^r_{A,K'})\ \longto\ \CM_k(M^r_{A,K}).
\end{equation}
The construction directly implies that
\addtocounter{Thm}{1}
\begin{equation}
\label{MF:Deg}
J_{g\kern1pt*}\circ J_g^* \ =\ [K:gK'g^{-1}]\cdot\id.
\end{equation}
Note also that with given $K$ and $K'$ the morphism $J_g$ and the maps $J_g^*$ and $J_{g\kern1pt*}$ depend only on the coset~$Kg$. Furthermore, an explicit calculation that we leave to the reader shows that
\addtocounter{Thm}{1}
\begin{equation}
\label{MF:PushComp}
J_{g\kern1pt*} \circ J_{g'*}  = J_{gg'*}
\end{equation}
whenever the formula makes sense. 

%%%%%%%%%%%%%%%%%%%%%%%

\medskip
Now consider an element $g\in\GL_r(\BA_F^f)$ with coefficients in $\hat A$ and an open compact subgroup $K\subset\GL_r(\hat A)$, bearing no particular relation with each other. Then with $K' := K \cap\nobreak g^{-1}Kg$ the pair of morphisms
\addtocounter{Thm}{1}
\begin{equation}
\label{HO:HC}
\xymatrix{
M^r_{A,K} & \ar[l]_{J_1} M^r_{A,K'} \ar[r]^{J_g} & M^r_{A,K} \\}
\end{equation}
is called the \emph{Hecke correspondence on $M^r_{A,K}$ associated to~$g$}. The composite map
\addtocounter{Thm}{1}
\begin{equation}
\label{HO:HO}
\xymatrix{
\llap{$T_g\!:\ \ $}\CM_k(M^r_{A,K}) \ar[r]^-{J^*_1} & 
\CM_k(M^r_{A,K'}) \ar[r]^-{J_{g\kern1pt*}} & \CM_k(M^r_{A,K}) \\}
\end{equation}
is called the \emph{Hecke operator on $\CM_k(M^r_{A,K})$ associated to~$g$}. It depends only on the double coset $KgK$. The composites of Hecke operators are calculated as follows:

\begin{Prop}
\label{HO:Comp}
For any $g$, $g'\in\GL_r(\BA_F^f)$ with coefficients in $\hat A$ and any open compact subgroup $K\subset\GL_r(\hat A)$ the Hecke operators on $\CM_k(M^r_{A,K})$ satisfy
$$T_{g'}\circ T_g \ =\ \sum_{g''} \; 
\bigl[ K\cap g^{\prime\prime-1}Kg'' : K\cap g^{-1}Kg\cap g^{\prime\prime-1}Kg'' \bigr]
\cdot T_{g''}$$
where $g''$ runs through a set of representatives of the double quotient
$$(K\cap g'Kg^{\prime-1}) \backslash g'Kg / (g^{-1}Kg\cap K).$$
\end{Prop}

\begin{Proof}
(Sketch) Consider the following diagram, where the disjoint union is extended over the set of representatives mentioned above:
$$\xymatrix@!C@C-60pt@R+10pt{
&& \bigsqcup_{g''} M^r_{A,\kern1pt K\cap\,g^{-1}\!K\kern-1pt g\,\cap\,g^{\prime\prime-1}\!K\kern-1pt g''} 
  \ar[dl]_{(J_1)} \ar[dr]^{(J_{g^{\prime-1}g''})} 
  \ar@/_35pt/[ddll]_{(J_1)}  \ar@/^35pt/[ddrr]^{(J_{g''})}
  && \\
& M^r_{A,\kern1pt K\,\cap\,g^{-1}\!K\kern-1pt g} \ar[dl]_{J_1} \ar[dr]^{J_g} && 
  M^r_{A,\kern1pt K\,\cap\,g^{\prime-1}\!K\kern-1pt g'} \ar[dl]_{J_1} \ar[dr]^{J_{g'}} & \\
M^r_{A,K} && M^r_{A,K} && M^r_{A,K} \\}$$
Since $Kg^{\prime-1}g'' = Kg$ for all $g''$ in question, the parallelogram in the middle commutes. We claim that it is cartesian over the open dense subset of $M^r_{A,K}$ over which both $J_g\!: M^r_{A,K\cap g^{-1}Kg} \to M^r_{A,K}$ and $J_1\!: M^r_{A,K\cap g^{\prime-1}Kg'} \to M^r_{A,K}$ are \'etale. This follows from a standard calculation using the analytic description (\ref{DMV:DQ}) and (\ref{MI:R1}) that we find too awkward to reproduce here. Since sections of invertible sheaves over integral schemes are determined by their restrictions to open dense subschemes, from the definitions of the maps we deduce that the composite
$$\xymatrix@C+57pt{
\CM_k\bigl(M^r_{A,K\cap g^{-1}Kg}\bigr) \ar[r]^-{J_{g\kern1pt*}} & 
\CM_k\bigl(M^r_{A,K}\bigr) \ar[r]^-{J_1^*} & 
\CM_k\bigl(M^r_{A,K\cap g^{\prime-1}Kg'}\bigr)\\}$$
along the lower edge of the parallelogram is equal to the sum of the composites
$$\xymatrix@C+20pt{
\CM_k\bigl(M^r_{A,K\cap g^{-1}Kg}\bigr) \ar[r]^-{J_1^*} & 
\CM_k\bigl(M^r_{A,\kern1pt K\cap\,g^{-1}\!K\kern-1pt g\,\cap\,g^{\prime\prime-1}\!K\kern-1pt g''}\bigr) 
\ar[r]^-{J_{g^{\prime-1}g''\kern1pt*}} & 
\CM_k\bigl(M^r_{A,K\cap g^{\prime-1}Kg'}\bigr)\\}$$
along the upper edge of the parallelogram. Thus $T_{g'}\circ T_g = J_{g'\kern1pt*} \circ J^*_1 \circ J_{g\kern1pt*} \circ J^*_1$ is the sum of the composites in the top row of this diagram:
$$\xymatrix@C-40pt{
\CM_k\bigl(M^r_{A,K}\bigr) \ar[rr]^-{J_1^*} && 
\CM_k\bigl(M^r_{A,\kern1pt K\cap\,g^{-1}\!K\kern-1pt g\,\cap\,g^{\prime\prime-1}\!K\kern-1pt g''}\bigr) 
\ar[rr]^-{J_{g''*}} && 
\CM_k\bigl(M^r_{A,K}\bigr) \\
& \CM_k\bigl(M^r_{A,\kern1pt K\cap\,g^{\prime\prime-1}\!K\kern-1pt g''}\bigr)
\ar[ur]^-{J_1^*} \ar@{<-}[ul]_{J_{1*}} 
% \ar[rr]
&& \CM_k\bigl(M^r_{A,\kern1pt K\cap\,g^{\prime\prime-1}\!K\kern-1pt g''}\bigr)
\ar@{<-}[ul]_-{J_1^*} \ar[ur]^{J_{g''*}} & \\}$$
%$$\xymatrix@C+20pt{
%\CM_k\bigl(M^r_{A,K}\bigr) \ar[r]^-{J_1^*} & 
%\CM_k\bigl(M^r_{A,\kern1pt K\cap\,g^{-1}\!K\kern-1pt g\,\cap\,g^{\prime\prime-1}\!K\kern-1pt g''}\bigr) 
%\ar[r]^-{J_{g''\kern1pt*}} & 
%\CM_k\bigl(M^r_{A,K}\bigr) \\
%& \CM_k\bigl(M^r_{A,\kern1pt K\cap\,g^{\prime\prime-1}\!K\kern-1pt g''}\bigr)
%\ar@<10pt>[u]^-{J_1^*} \ar@<-10pt>@{<-}[u]_{J_{1\kern1pt*}} 
%\ar@{<-}[ul]^-{J_1^*} \ar[ur]_{J_{g''\kern1pt*}} & \\}$$
The indicated factorizations and (\ref{MF:Deg}) now yield the desired formula.
\end{Proof}

%%%%%%%%%%%%%%%%%%%%%%%%%%%%%%%%%%%%%%%%%%%%%%%%%%%%%%%%%%%%%%%%%%%%%%%%%%%%%%%%%%%%%%%%%%%%%

\section{The special case $\BF_q[t]$ and level $(t)$}
\label{Fqtt}

%%%%%%%%%%%%%%%%%%%%%%%

Throughout this section we consider the case where $A:=\BF_q[t]$ for a finite field $\BF_q$ of cardinality $q=p^m$, the level is $K(t)$, and the rank $r\ge1$ is arbitrary. In this case the Satake compactification was already described by Kapranov \cite{Kapranov}; but we will obtain finer information about it. We set $V_r := \BF_q^r$ and identify it with the $\BF_q$-vector space $(t^{-1}A/A)^r$. 
% For any $\BF_q$-vector space $V$ we set $\circV \ :=\ V \setminus \{0\}$.

\medskip
Let $S$ be a scheme over $F=\BF_q(t)$, and let $(E,\phi,\lambda)$ be a Drinfeld $A$-module of rank $r$ with a level $(t)$ structure over~$S$. 
% Then $\phi$ is determined by the single element~$\phi_t$, which must be of the form $\sum_{i=0}^r \phi_{t,i} \tau^{im}$ with $\phi_{t,i}\in\Gamma(S,E^{1-q^i})$ and $\phi_{t,r}$ invertible.
By definition $\lambda$ is an $A$-linear isomorphism $\lambda\!:\ \underline{V_r} \stackrel\sim\to \Ker(\phi_t)\subset E$. We can view it equivalently as an $\BF_q$-linear map $V_r\to E(S)$ satisfying certain additional conditions. In particular it is fiberwise injective, i.e., the composite map $V_r\to E(S) \to E(s)$ is injective for every point $s\in S$. It turns out that $(E,\lambda)$ determines $\phi$ completely, because:

\begin{Prop}
\label{Fqtt:Recon}
For any line bundle $E$ over $S$ and any fiberwise injective $\BF_q$-linear map $\lambda\!: V_r\to E(S)$ there exists a unique homomorphism $\phi\!: A \to \End(E)$ turning $(E,\phi,\lambda)$ into a Drinfeld $A$-module of rank $r$ with level $(t)$ structure over~$S$. 
\end{Prop}

\begin{Proof}
The assertion is local on $S$, so we may assume that $E=\BG_{a,S}$ and $S=\Spec R$. Then $\lambda$ is an $\BF_q$-linear homomorphism $V_r\to R$ such that $\lambda(v)$ is invertible for all non-zero~$v$. Giving a Drinfeld $A$-module $\phi\!: A \to \End(\BG_{a,S})$ of rank $r$ is equivalent to giving the single element $\phi_t\in \End(\BG_{a,S}) = R[\tau]$, which must be of the form $\phi_t\ =\ \sum_{i=0}^r c_i\cdot\tau^{im}$ with $c_i\in R$ and $c_0=t$ and invertible highest coefficient~$c_r$. 

The level structure requires in addition that $\Ker(\phi_t) = \lambda\bigl(\underline{V_r}\bigr)$. We claim that this is equivalent to
\addtocounter{Thm}{1}
\begin{equation}
\label{Fqtt:Phit}
\phi_t(X)\ =\ 
t\cdot X\cdot \!\!\!\!\prod_{v\in V_r\setminus\{0\}} \Bigl( 1 - \frac{X}{\lambda(v)} \Bigr).
\end{equation}
Indeed, the right hand side vanishes to first order at all prescribed zeros of~$\phi_t$; hence $\phi_t(X)$ must be a multiple of the right hand side, say by the element $f(X) \in R[X]$. Since both sides of (\ref{Fqtt:Phit}) are polynomials of degree $q^r$ in~$X$ and possess invertible highest coefficients, this $f$ must be constant. As the coefficient of $X$ on both sides is the unit~$t$, we must in fact have $f=1$. This shows that the equality (\ref{Fqtt:Phit}) is necessary. It is also clearly sufficient.

It remains to show that (\ref{Fqtt:Phit}) actually does define a Drinfeld $A$-module of rank $r$ with level structure~$\lambda$. For this write the right hand side as a unit times $\prod_{v\in V_r} (X-\nobreak\lambda(v))$. Since $\lambda$ is $\BF_q$-linear, it is well-known \cite[Cor.$\,$1.2.2]{GossBS} that any such polynomial is $\BF_q$-linear; therefore $\phi_t(X)=\sum_{i=0}^r c_i\cdot X^{q^i}$ with $c_i\in R$. The formula (\ref{Fqtt:Phit}) also shows that $c_0=t$ and that $c_r$ is invertible. Thus $\phi_t(X)$ extends to a unique Drinfeld $A$-module $\phi\!: A \to R[\tau]$ of rank~$r$. By the preceding remarks $\lambda$ already defines an $\BF_q$-linear isomorphism $\underline{V_r} \stackrel\sim\longto \Ker(\phi_t)$. Since $t$ annihilates both sides, the isomorphism is then in fact $A$-linear, as desired.
\end{Proof}

%%%%%%%%%%%%%%%%%%%%%%%
\medskip

Recall that the projective space $\BP^{r-1}$ represents the functor that to any scheme $S$ associates the set of isomorphism classes of tuples $(E,e_1,\ldots,e_r)$ consisting of a line bundle $E$ on $S$ and sections $e_i\in E(S)$ which generate~$E$. Let $\Omega_r$ denote the open subvariety of $\BP^{r-1}_{\BF_q}$ obtained by removing all $\BF_q$-rational hyperplanes. (This definition is entirely analogous to the definition of the non-archimedean Drinfeld period domain $\Omega^r$ associated to the local field~$F_\infty$. We hope that the confusion be limited by the fact that the new $\Omega_r$ is not used outside the present section.) Then $\Omega_r$ represents the functor that to any scheme $S$ over $\BF_q$ associates the set of isomorphism classes of tuples $(E,e_1,\ldots,e_r)$ consisting of a line bundle $E$ on $S$ and sections $e_i\in E(S)$ which are fiberwise $\BF_q$-linearly independent. Giving such sections $e_i$ is equivalent to giving the $\BF_q$-linear map $\lambda\!: V_r\to\nobreak E(S)$, $(\alpha_1,\ldots,\alpha_r) \mapsto \sum_i \alpha_i e_i$, which must be fiberwise injective. Thus the pullback $\Omega_{r,F}$ of $\Omega_r$ to $\Spec F$ represents the functor that to any scheme $S$ over $F$ associates the set of isomorphism classes of pairs $(E,\lambda)$ consisting of a line bundle $E$ on $S$ and fiberwise injective $\BF_q$-linear map $\lambda\!: V_r\to E(S)$. In view of Proposition \ref{Fqtt:Recon} we obtain an isomorphism for the moduli space of Drinfeld modules
\addtocounter{Thm}{1}
\begin{equation}
\label{Fqtt:ModOmega}
M^r_{\BF_q[t],K(t)}\ \cong\  \Omega_{r,F}.
\end{equation}
Somewhat surprisingly (compare \cite{GekelerSat}), its Satake compactification is \emph{not} $\BP^{r-1}_F$ when $r\ge3$. That was the motivation for the article \cite{PinkSchieder}, where another compactification of $\Omega_r$ was constructed and studied in detail. This compactification is defined as follows.

\medskip
Let $S_r$ denote the symmetric algebra of $V_r$ over $\BF_q$, which is a polynomial ring in $r$ independent variables. Let $K_r$ denote its field of quotients, and let $R_r \subset K_r$ be the $\BF_q$-subalgebra generated by the elements $\frac{1}{v}$ for all $v\in V_r\setminus\{0\}$. Turn $R_r$ into a graded $\BF_q$-algebra by declaring each $\frac{1}{v}$ to be homogeneous of degree~$1$. Let $RS_r\subset K_r$ denote the subalgebra generated by $R_r$ and~$S_r$.

Then $\BP^{r-1}_{\BF_q} = \Proj S_r$, and the localization $RS_r$ of $S_r$ corresponds to the open subscheme~$\Omega_r$. By construction $RS_r$ is also a localization of $R_r$; and so $\Omega_r$ is also an open subscheme of the projective scheme $Q_r := \Proj R_r$. This is the variety that we are interested~in. With (\ref{Fqtt:ModOmega}) it follows that $Q_{r,F}$ is a projective compactification of $M^r_{\BF_q[t],K(t)}$. Let $\CO(1)$ denote the natural very ample invertible sheaf on it, whose space of global sections contains the elements $\frac{1}{v}\in R_r$.

\begin{Thm}
\label{Fqtt:Sat}
The variety $Q_{r,F}$ is a Satake compactification of $M^r_{\BF_q[t],K(t)}$. 
The dual of the relative Lie algebra of its universal family is $\CO(1)$.
\end{Thm}

\begin{Proof}
{}From \cite[Thm.$\,$1.9]{PinkSchieder} we know that $Q_{r,F}$ is a normal integral algebraic variety. Let $\bar E$ denote the line bundle on it whose sheaf of sections is $\CO(-1)$. By construction the dual of its relative Lie algebra is~$\CO(1)$. By Definition \ref{SC:SatDef} it remains to construct a homomorphism $\bar\phi\!: A\to\End(\bar E)$ which turns $(\bar{E},\bar{\phi})$ into a weakly separating generalized Drinfeld $A$-module over $Q_{r,F}$ whose restriction to $M^r_{\BF_q[t],K(t)}$ is the given universal family $(E,\phi)$.

Since $\frac{1}{v}$ is a section of $\CO(1)$, for any local section $f$ of $\bar E$ the product $\frac{1}{v}\cdot f$ is a local section of the structure sheaf~$\CO_{Q_{r,F}}$. Thus if we plug $X=f$ into the polynomial 
\addtocounter{Thm}{1}
\begin{equation}
\label{Fqtt:Phit2}
\bar\phi_t(X)\ :=\ t\cdot X\cdot \!\!\!\!\prod_{v\in V_r\setminus\{0\}} 
\Bigl( 1 - \frac{1}{v}\cdot X \Bigr),
\end{equation}
all but the single factor $X$ turn into sections of $\CO_{Q_{r,F}}$, and so the result is again a local section of~$\bar E$. Therefore $\bar\phi_t$ defines a morphism of algebraic varieties $\bar E\to\bar E$ over~$Q_{r,F}$. Restricted to $\Omega_{r,F}\subset Q_{r,F}$, the sections $\frac{1}{v}$ become invertible and their inverses $v$ are precisely the non-zero elements of the $\BF_q$-subspace $V_r \subset \bar E(\Omega_{r,F})$. Thus comparison with (\ref{Fqtt:Phit}) shows that the restriction $(\bar E,\bar\phi_t)|\Omega_{r,F}$ is isomorphic to $(E,\phi_t)$. Since $Q_{r,F}$ is integral, it follows that $\bar\phi_t$ is $\BF_q$-linear of degree $\le rm$ as a (non-commutative) polynomial in~$\tau$ everywhere. As before it extends to a unique $\BF_q$-algebra homomorphism $\bar\phi\!: A\to\End(\bar E)$ whose restriction to $\BF_q$ is induced by the given embedding $\BF_q\into F$, and by construction we have $(\bar E,\bar\phi)|\Omega_{r,F} \cong (E,\phi)$.

% By construction $\bar\phi_t(X)$ has degree $|V_r|=q^r=p^{rm}$ in $X$; hence $\bar\phi_t$ has degree $rm$ in~$\tau$. This shows 
Thus the condition (c) in Definition \ref{GDM:GenDrinDef} has been shown, and (a) holds because the coefficient of $X$ in $\bar\phi_t(X)$ is~$t$. For (b) recall that the elements $\frac{1}{v}$ generate the ring~$R_r$, and so the corresponding sections generate the sheaf $\CO(1)$. Thus at every point on $Q_{r,F}$, at least one of these sections is non-zero, and so the polynomial $\bar\phi_t(X)$ is not just linear in~$X$. This implies~(b); hence $(\bar E,\bar\phi)$ is a generalized Drinfeld $A$-module of rank $\le r$.

Finally, consider a set of points of $Q_{r,F}$ over a field~$L$ at which the fibers of $(\bar E,\bar\phi)$ are all isomorphic. Then the zero sets of the polynomial $\bar\phi_t(X)$ at these points are equal up to multiplication by an element of~$L^\times$. By the definition of $\bar\phi_t(X)$ this means that the values of the generators $\frac{1}{v}$ at these points are equal up to a permutation and joint multiplication by an element of~$L^\times$. As the number of possible permutations is finite, so is the number of points; hence $(\bar E,\bar\phi)$ is weakly separating, as desired.
\end{Proof}

%%%%%%%%%%%%%%%%%%%%%%%

\medskip
In the rest of this section we abbreviate
\begin{eqnarray*}
M^r &:=& M^r_{\BF_q[t],K(t)}     \ \ \cong\ \ \Omega_{r,F},\\
\OM^r &:=& \OM^r_{\BF_q[t],K(t)} \ \ \cong\ \      Q_{r,F}.
\end{eqnarray*}

\begin{Rem}
\label{Fqtt:SatRem}
\rm Section 6 of \cite{PinkSchieder} gives a modular interpretation of~$Q_r$, which in the context of the present paper can be viewed as describing generalized Drinfeld $\BF_q[t]$-modules of rank $\le r$ with a level $(t)$ structure. It would be interesting to know whether the Satake compactification possesses a similar modular interpretation in general.
\end{Rem}

\begin{Rem}
\label{Fqtt:StratRem}
\rm Section 7 of \cite{PinkSchieder} describes a stratification of $Q_r$ whose strata are isomorphic to $\Omega_{r'}$ and indexed by $\BF_q$-subspaces $V'\subset V_r$ of dimension~$r'$  for all $1\le r'\le r$. This yields a stratification of $\OM^r$ by copies of Drinfeld modular varieties~$M^{r'}$, which was also described by Kapranov \cite[Thm.$\,$1.1]{Kapranov}.
\end{Rem}

\begin{Thm}
\label{Fqtt:SatSing}
(\cite[Thm.$\,$7.4]{PinkSchieder})
The singular locus of $\OM^r$ is the union of all strata of codimension $\ge\nobreak2$. In particular $\OM^r$ is singular whenever $r\ge3$.
\end{Thm}

\begin{Rem}
\label{Resol}
\rm Section 9 of \cite{PinkSchieder} constructs a resolution of singularities of $\OM^r$. 
\end{Rem}
% Next we turn to modular forms on~$M^r$.

\begin{Thm}
\label{Fqtt:ModRV}
The ring of modular forms $\CR(M^r)$ is isomorphic to $R_r\otimes_{\BF_q}F$ as a graded $F$-algebra. Moreover it is Cohen-Macaulay, and so $\OM^r$ is Cohen-Macaulay.
\end{Thm}

\begin{Proof}
By Definition \ref{MF:MFDef} and Theorem \ref{Fqtt:Sat} we have $\CM_k(M^r) = H^0\bigl(Q_{r,F},\CO(k)\bigr)$, and by \cite[Cor.$\,$4.4]{PinkSchieder} the latter is the degree $k$ part of $R_r\otimes_{\BF_q}F$. This proves the first assertion. The rest is the content of \cite[Thms.$\,$1.7,$\,$1.9]{PinkSchieder}.
\end{Proof}

\begin{Thm}
\label{Fqtt:Dim}
For all $k\ge0$ we have 
$$\dim_F \CM_k(M^r) \ =
\sum_{i_1, \ldots, i_{r-1} \in \{0,1\}} q^{\sum_{\nu} \nu \cdot i_{\nu}} \cdot \binom{k}{\sum_{\nu} i_{\nu}}.$$
\end{Thm}

\begin{Proof}
Theorem \ref{Fqtt:ModRV} and \cite[Thm.$\,$1.8]{PinkSchieder}.
\end{Proof}

%%%%%%%%%%%%%%%%%%%%%%%

\medskip
Finally, let $\omega$ denote the canonical sheaf on the regular locus $\OM^{r,\rm reg}$ of~$\OM^r$. Let $\partial M^r := \OM^r\setminus M^r$ be the boundary with the unique structure as reduced closed subscheme. Then by \cite[Thm.$\,$5.1]{PinkSchieder} there is an isomorphism
\addtocounter{Thm}{1}
\begin{equation}
\label{Fqtt:Can}
\CL_{\OM^r}^r|\OM^{r,\rm reg}\ \cong\ \omega(2\cdot\partial M^r).
\end{equation}
Since $\OM^r$ is normal, every section of $\CL_{\OM^r}^r$ over $\OM^{r,\rm reg}$ extends to $\OM^r$. Thus the space of modular forms of weight $r$ in this case is isomorphic to the space of top differentials on $\OM^{\rm reg}$ with at most double poles at the boundary:

\begin{Thm}
\label{Fqtt:Diff}
There is a natural isomorphism
$$\CM_r(M^r)\ \cong\ H^0(\OM^{r,\rm reg},\omega(2\cdot\partial M^r)).$$
\end{Thm}

For all the results from Remark \ref{Fqtt:SatRem} through Theorem \ref{Fqtt:Diff}, it would be natural to ask whether they generalize to arbitrary $K$ and arbitrary~$A$. A slight generalization of the dimension formula from Theorem \ref{Fqtt:Dim} is given in Theorem \ref{Fqt+:Dim} below.

%%%%%%%%%%%%%%%%%%%%%%%%%%%%%%%%%%%%%%%%%%%%%%%%%%%%%%%%%%%%%%%%%%%%%%%%%%%%%%%%%%%%%%%%%%%%%

\section{The special case $\BF_q[t]$ and level containing $(t)$}
\label{Fqt+}

In this section we keep $A:=\BF_q[t]$ and consider certain subgroups $K\subset\GL_r(\hat A)$ containing $K(t)$. In view of Propositions \ref{SC:SatQuot} and \ref{MF:Quot} the results will be obtained from those of the preceding section by taking quotients and invariants. Set
\begin{eqnarray*}
K_1(t) &:=& \left\{ g\in\GL_r(\hat A) \;\left|\; g\equiv 
\begin{pmatrix} 1 && * \\[-6pt] & \!\!\!\ddots\!\!\! & \\[-6pt] 0 && 1 \end{pmatrix}
\mod(t) \right.\right\}, \\[5pt]
K'(1) &:=& \bigl\{ g\in\GL_r(\hat A) \bigm| \det(g)\equiv1 \mod(t) \bigr\}, \\[5pt]
K(1) &:=& \GL_r(\hat A).
\end{eqnarray*}
{}From Theorem \ref{Fqtt:Sat} and Proposition \ref{SC:SatQuot} and \cite[Thm.$\,$3.2, Prop.$\,$3.3]{PinkSchieder} we deduce:

\begin{Thm}
\label{Fqt+:Sat}
\begin{itemize}
\item[(a)] 
$\OM^r_{\BF_q[t],K_1(t)} \cong \BP^{r-1}_F$.
\item[(b)] 
$\OM^r_{\BF_q[t],K'(1)}$ is a weighted projective space of weights $q{-}1,\ldots,\,q^{r-1}{-}1,\,\smash{\frac{q^r-1}{q-1}}$.
\item[(c)] 
$\OM^r_{\BF_q[t],K(1)}$ is a weighted projective space of weights $q{-}1,\,q^2{-}1,\ldots,\,q^r{-}1$.
\end{itemize}
\bigskip\noindent
In (a) the Satake compactification is smooth, in the other cases it is singular if $r\ge3$.
\end{Thm}

{}From Theorem \ref{Fqtt:ModRV} and Proposition \ref{MF:Quot} and \cite[Thm.$\,$3.1]{PinkSchieder} we deduce:

\begin{Thm}
\label{Fqt+:MF}
The following rings are generated over $F$ by $r$ algebraically independent modular forms of the indicated weights:
\begin{itemize}
\item[(a)] $\CR\bigl(M^r_{\BF_q[t],K_1(t)}\bigr)$ with all weights $1$.
\item[(b)] $\CR\bigl(M^r_{\BF_q[t],K'(1)}\bigr)$ 
with weights $q{-}1,\ldots,\,q^{r-1}{-}1,\,\smash{\frac{q^r-1}{q-1}}$.
\item[(c)] $\CR\bigl(M^r_{\BF_q[t],K(1)}\bigr)$ 
with weights $q{-}1,\,q^2{-}1,\ldots,\,q^r{-}1$.
\end{itemize}
\end{Thm}

\begin{Rem}
\label{Fqt+:Breuer}
\rm More precisely, the ring in (c) is isomorphic to $R_r^{\GL_r(\BF_q)} \otimes_{\BF_q}F$, and the proof of \cite[Thm.$\,$3.1]{PinkSchieder} shows that its generators
% of $\CR\bigl(M^r_{\BF_q[t],K(1)}\bigr)\cong R_r^{\GL_r(\BF_q)} \otimes_{\BF_q}F$ 
correspond to the coefficients of $X^{q^i}$ for $1\le i\le r$ in the polynomial $\bar\phi_t(X)$ from (\ref{Fqtt:Phit2}). In other words, the coefficients of the universal Drinfeld $\BF_q[t]$-module of rank $r$ form algebraically independent generators of the ring of modular forms on $M^r_{\BF_q[t],K(1)}$. 
This was pointed out to me by Breuer, and seems to be implicit in Kapranov \cite[Rem.$\,$1.6]{Kapranov}.
% , an observation due to Breuer.
\end{Rem}

Finally, from Theorem \ref{Fqtt:ModRV} and Proposition \ref{MF:Quot} and \cite[Thm.$\,$4.1]{PinkSchieder} we deduce:

\begin{Thm}
\label{Fqt+:Dim}
For any $K$ satisfying $K(t)\subset K\subset K_1(t)$ and all $k\ge0$ we have 
$$\dim_F \CM_k\bigl(M^r_{\BF_q[t],K}\bigr)\ =\ 
  \sum_{s=1}^r\ \frac{\,|K\backslash\GL_r(\hat A)/J_s(\hat A)|\,}{\prod_{i=1}^s(q^i-1)} \cdot \binom{k-1}{s-1},$$
where $J_s\subset\GL_r$ is the subgroup of elements whose first $s$ columns coincide with those of the identity matrix.
\end{Thm}

It is natural to ask whether the same dimension formula holds for any fine open compact subgroup~$K$ and whether it generalizes to arbitrary~$A$.
% $K\subset K_1(t)$.

%%%%%%%%%%%%%%%%%%%%%%%%%%%%%%%%%%%%%%%%%%%%%%%%%%%%%%%%%%%%%%%%%%%%%%%%%%%%%%%%%%%%%%%%%%%%%

%%%%%%%%%%%%%%%%%%%%%%%%%%%%%%%%%%%%%%%%%%%%%%%%%%%%%%%%%%%%%%%%%%%%%%%%%%%%%%%%%%%%%%%%%%%%%

%%%%%%%%%%%%%%%%%%%%%%%%%%%%%%%%%%%%%%%%%%%%%%%%%%%%%%%%%%%%%%%%%%%%%%%%%%%%%%%%%%%%%%%%%%%%%
%
%\begin{center}
%\rule{8cm}{0.01cm}
%\end{center}
%
%\begin{minipage}[t]{8cm}{\small
%Department of Mathematical Sciences \\
%University of Stellenbosch \\
%Stellenbosch, 7600 \\
%South Africa \\
%fbreuer@sun.ac.za}
%\end{minipage}
%%
%%\hspace{2cm}
%\begin{minipage}[t]{8cm}{\small
%Department of Mathematics \\
%ETH Z\"urich\\
%8092 Z\"urich\\
%Switzerland \\
%pink@math.ethz.ch}
%\end{minipage}

\end{document}